\documentclass[times]{MyArticle}

\usepackage[colorlinks,bookmarksopen,bookmarksnumbered,citecolor=red,urlcolor=red]{hyperref}

\begin{document}

\runningheads{P.~Suchde \emph{et~al.}}{Flux Conserving Meshfree Method}

\title{A Flux Conserving Meshfree Method for Conservation Laws}

\author{Pratik Suchde \affil{1,2}\corrauth, %
J\"org Kuhnert \affil{1}, %
Simon Schr\"oder \affil{1}and
Axel Klar \affil{2} %
}

\address{\affilnum{1}Fraunhofer ITWM, 67663 Kaiserslautern, Germany \break %
\affilnum{2}Department of Mathematics, University of Kaiserslautern, 67663 Kaiserslautern, Germany %
}

\corraddr{P. Suchde. E-mail: pratik.suchde@itwm.fraunhofer.de}

\begin{abstract}%
Lack of conservation has been the biggest drawback in meshfree generalized finite difference methods~(GFDMs). In this paper, we present a novel modification of classical meshfree GFDMs to include local balances which produce an approximate conservation of numerical fluxes. This numerical flux conservation is done within the usual moving least squares framework. Unlike Finite Volume Methods, it is based on locally defined control cells, rather than a globally defined mesh. We present the application of this method to an advection diffusion equation and the incompressible Navier--Stokes equations. Our simulations show that the introduction of flux conservation significantly reduces the errors in conservation in meshfree GFDMs.
\end{abstract}


\maketitle

\section{Introduction}
Generation and management of meshes is often the most difficult and time consuming part of numerical simulation procedures. This is further compounded for complex, time-dependent geometries. The efficiency of mesh generation determines the overall accuracy and robustness of the simulation process. To avoid the task of meshing, several classes of meshless or meshfree methods have been developed. Meshfree methods use the numerical basis of a set of nodes, which need not be regularly distributed, to cover the computational domain. These nodes could either be mass-carrying particles or numerical points. For each node, the only connectivity information required is a local set of neighbouring nodes over which approximations are carried out. 

Smoothed Particle Hydrodynamics~(SPH) \cite{Liu2003SPH} is one of the first and most widely used meshless methods. In SPH, the computational domain is discretized by particles which carry mass. Thus, mass conservation is directly guaranteed. One of the biggest drawbacks of SPH is the difficulty in enforcing boundary conditions \cite{Li2002SPH, Liu2011SPH, Randles1996SPH}. While a lot of work has been done to address this issue, the SPH formulation does not naturally include treatment of most boundary conditions and extra effort is needed to enforce them. One such method is the addition of ghost or dummy particles outside the computational domain \cite{Li2002SPH}. Another major drawback of SPH in its classical formulation is the so-called particle inconsistency issue which results in the absence of a valid approximation order. The original SPH does not even have $C^0$ particle consistency for irregularly distributed particles and boundary particles. Attempts to solve this problem include the kernel renormalization \cite{Chen2000SPH} which comes at the price of the loss of certain conservation properties for momentum and energy. In addressing these issues, the Finite Pointset Method~(FPM) evolved out of SPH \cite{KuhnertThesis}. FPM \cite{Jefferies2015, Kuhnert2014, Tiwari2016, Tiwari2002} is a meshfree method based on moving least squares~(MLS) approximations such that boundary conditions can be directly enforced, without any extra effort. While SPH uses a set of mass-carrying particles as a numerical basis, FPM and other meshfree GFDMs (for example, \cite{Katz2010,Liszka1980,Praveen2007,Sridar2003}) use a set of numerical points, referred to as a point cloud. These numerical points do not carry mass, and thus, points can easily be added or removed during a simulation \cite{Jefferies2015}, which is not possible in the particle-based SPH. Since point clouds can easily be modified locally, FPM is more adaptive than SPH, in terms of the discretization of the domain. However, these advantages come with the drawback of a lack of formal conservation. Despite the lack of conservation, the FPM has been shown to be a robust method with many practical applications \cite{Drumm2008,Jefferies2015, Tiwari2016, Uhlmann2013}. The FPM is also used as the numerical basis of two commercially used meshfree simulation tools: NOGRID \cite{Moller2007} and the meshfree module of VPS-PAMCRASH \cite{Tramecon2013}. It must be noted that in the meshfree context, the acronym FPM is often a confusing one. It is used to represent not only the aforementioned Finite Pointset Method, but also the Finite Particle Method \cite{LiuFPM} and the well-known Finite Point Method \cite{OnateFPM}. Thus, we henceforth drop the acronym FPM and refer to the Finite Pointset Method under the umbrella term of GFDM. 

Despite the name, most meshfree methods are not completely devoid of meshes. While meshfree methods lack the use of a \emph{predefined} mesh for domain discretization, many meshfree methods use an easily generable mesh such that the solution is not too heavily dependent on the quality of that mesh~(see \cite{Liu2009Meshfree} for details). Several methods, which solve partial differential equations in their weak formulation, use a so-called background mesh for numerical integration of the weak-forms. The element-free Galerkin method~(EFG) \cite{Belytschko1994} uses a global background mesh, while the meshless local Petrov-Galerkin method~(MLPG) \cite{Atluri1998} requires only local background cells. Meshfree particle methods, such as SPH, use a predefinition of mass particles, which often requires some kind of a mesh. One of the methods to solve the problem of particle distortion in SPH is by regularization, which requires a re-meshing \cite{Tu2012}. Several methods in the family of point interpolation methods~(PIM) \cite{Liu2013PIM} use tessellations or background cells in so-called T-schemes for the selection of support domains. Point-based meshfree GFDMs are often referred to as ``truly meshfree". However, they too use some notion of a mesh. They often use local triangulations for accurate post-processing operations or for the identification of points on the free surface of fluid flow simulations. In each of these cases, the meshes used are a lot easier to generate and do not impose the restrictions that meshes in classical meshed methods do.

A few attempts have been made to overcome the lack of conservation in meshfree GFDMs \cite{Chiu2012, Chiu2011}, but they lose the completely local nature of the differential operators. They introduce a symmetry of differential operator coefficients which leads to a discrete divergence theorem. However, that results in a need to solve for differential operators in a global system. That increases simulation time significantly and makes those methods extremely inefficient on finely discretized point clouds and moving point clouds which require a re-computation of numerical differential operators at every time step.

Finite Volume Methods~(FVMs) have the feature of local conservation of numerical fluxes from one control volume to its neighbouring one \cite{Eymard2000}. In this paper, we aim to use local control cells, in a manner similar to FVMs, in a meshfree framework to introduce an approximate conservation of numerical fluxes. A local control cell at a point is formed by the Delaunay tessellation of points in its support domain. The differential operators are computed using an MLS approach which guarantees usual monomial consistency coupled with a conservation of numerical fluxes for specific fields. For modeling fluid flow, we use a fully Lagrangian framework. Since point clouds can be modified locally, and the control cells are also defined locally, this method does not face the problem of cell distortion present in Lagrangian moving mesh methods. It must be noted that the idea of generalizing control volumes used by FVMs in the context of meshfree methods is not a novel one. The Finite-Volume Particle Method~(FVPM) \cite{Hietel2005,Hietel2000} generalizes control volumes to include those which need not be disjoint. In FVPM, grid generation is replaced by an expensive integration of partition of unity functions. The Meshless-Finite Volume Method~(MFVM) \cite{Hopkins2015} generalizes Voronoi-based moving mesh methods. MFVM uses a volume partitioning that amounts to a Voronoi tessellation with edges smoothed. Unlike both the FVPM and MFVM, the method presented here is a strong-form method.

The paper is organized as follows. In Section~\ref{sec:FPM}, we give a brief overview of GFDMs. In Section~\ref{sec:Conservation}, we introduce our proposed method by presenting the details of control cells in a meshfree context and the construction of flux conserving differential operators. In Section~\ref{sec:AdvDiff} and \ref{sec:INSE}, we present the application of this method to an advection-diffusion equation and the incompressible Navier--Stokes equations respectively. Finally, we conclude the paper in Section~\ref{sec:Conclusion} with a discussion of the proposed method.
%
%
%
%
%
%
\section{Meshfree Generalized Finite Difference Methods}
\label{sec:FPM}
For the Finite Pointset Method and other similar meshfree GFDMs, the computational domain $\Omega$, with boundary $\partial\Omega$, is discretized using a cloud of $N$ numerical points with positions $\vec{x}_i$, $i = 1, \dots, N$, which includes points both in the interior and on the boundary of the domain. The points are usually irregularly spaced. Each numerical point carries the necessary numerical data of the problem. Each point $i$ has a set of neighbouring points $S_i$ which contains $n_i$ points, including itself. The neighbourhood or support $S_i$ is determined as $S_i=\{\vec{x}_j : \|\vec{x}_j-\vec{x}_i\| \le \beta h\}$, where $h=h(\vec{x},t)$ is the radius of the support, referred to as smoothing length or interaction radius; and $\beta \leq 1$ is a positive constant. A numerical point is not a mass carrying particle and thus, points can easily be added or deleted to adapt the point cloud locally. The spatial distribution of points is described by three parameters: $h$, $r_{min}$ and $r_{max}$. It is ensured that no two points are closer than $r_{min}h$ and there exists at least one point in every possible sphere of radius $r_{max}h$ in the computational domain. Thus, the smoothing length $h$ also determines the spatial discretization size. $r_{min}$ and $r_{max}$ usually have values of approximately $0.2$ and $0.45$ respectively. Further details about point cloud organization and management, including the setup of the initial point cloud, can be found in \cite{Drumm2008,Jefferies2015}.

As the name suggests, numerical derivatives are approximated with a generalized finite difference approach. For a function $f$ defined at each numerical point $i=1,2,\dots,N$, its derivatives are approximated as
\begin{equation}
	\label{Eq:FPM_DiffOp}
	\partial^*_i f(\vec x)\approx \tilde{\partial}^*_if = \sum_{j\in S_i}c_{ij}^*f_j\,,
\end{equation}
where ${}^*=x,y, xx, \Delta, etc.$ represents the differential operator being approximated, $\partial^*_i$ represents the continuous $^*$-derivative at point $i$, and $\tilde{\partial}^*_i$ represents the discrete derivative. For a point $i$, the stencil coefficients $c_{ij}$ are found using a weighted least squares approach. The weighted sum of the stencil coefficients is minimized such that the derivatives of monomials $m\in\mathcal{M}$ up to a certain order, usually $2$, are exactly reproduced.
\begin{align}
	\sum_{j\in S_i}c_{ij}^{*}m_j &= \partial^*_i m \qquad \forall m\in\mathcal{M}\,,\label{Eq:Consistency}\\
	\text{min } J &= \sum_{j\in S_i} W_{ij}(c_{ij}^{*})^2\,, \label{Eq:BasicMin}
\end{align}
where $W$ is a weighting function. Throughout this paper, we use a Gaussian weighting function
\begin{equation}
	W_{ij}=\exp(-4 \frac{\|\vec{x}_j-\vec{x}_i\|^2}{h_i^2 + h_j^2})\,.
\end{equation}
For each point $i$, Eq.\,\eqref{Eq:Consistency} and Eq.\,\eqref{Eq:BasicMin} are solved to obtain the numerical differential operator, which can then be substituted in the partial differential equations being solved to obtain the spatial discretization. Henceforth, the differential operators found by Eq.\,\eqref{Eq:Consistency} and Eq.\,\eqref{Eq:BasicMin} will be referred to as classical GFDM.
\section{Conservation}
\label{sec:Conservation}
Consider a conservation law 
\begin{equation}
	\label{Eq:ConservationLaw}
	\frac{\partial\phi}{\partial t}+\nabla\cdot\mathbf{J}=0\,,
\end{equation}
with $\mathbf{J}=\mathbf{J}(\phi)$. For sufficiently smooth $\phi$ and $\mathbf{J}$, integrating Eq.\,\eqref{Eq:ConservationLaw} over the entire domain, and an application of the divergence theorem leads to the integral form of the conservation law,
\begin{equation}
	\label{Eq:EnergyConservation}
	\frac{d}{dt}\int_\Omega\phi\, dV=-\int_{\partial\Omega} \vec{n}\cdot \mathbf{J}\,dA \,.
\end{equation}
Physically, Eq.\,\eqref{Eq:EnergyConservation} can be stated as the rate of change of energy within the domain should be the same as the energy flux across its boundary, when no energy source or sink is present. Finite volume methods \cite 
{Eymard2000} use a local balance on each discretization cell or control volume. FVMs solve Eq.\,\eqref{Eq:ConservationLaw} by directly enforcing a local, discrete version of Eq.\,\eqref{Eq:EnergyConservation}.

Meshfree GFDMs, on the other hand, usually directly solve the strong form Eq.\,\eqref{Eq:ConservationLaw}. The lack of a discrete divergence theorem causes the absence of a discrete form of Eq.\,\eqref{Eq:EnergyConservation}. This combined with the local nature of the differential operators, Eq.\,\eqref{Eq:FPM_DiffOp}, results in conservation not being ensured at a discrete level. We now present a novel modification of GFDM differential operators to enforce a local discrete divergence theorem. This results in an approximate discrete form of global conservation, Eq.\,\eqref{Eq:EnergyConservation}.
\subsection{A Meshfree Control Cell}
\label{sec:ControlCell}
\begin{figure}
  \centering
  \includegraphics[width=0.4\textwidth]{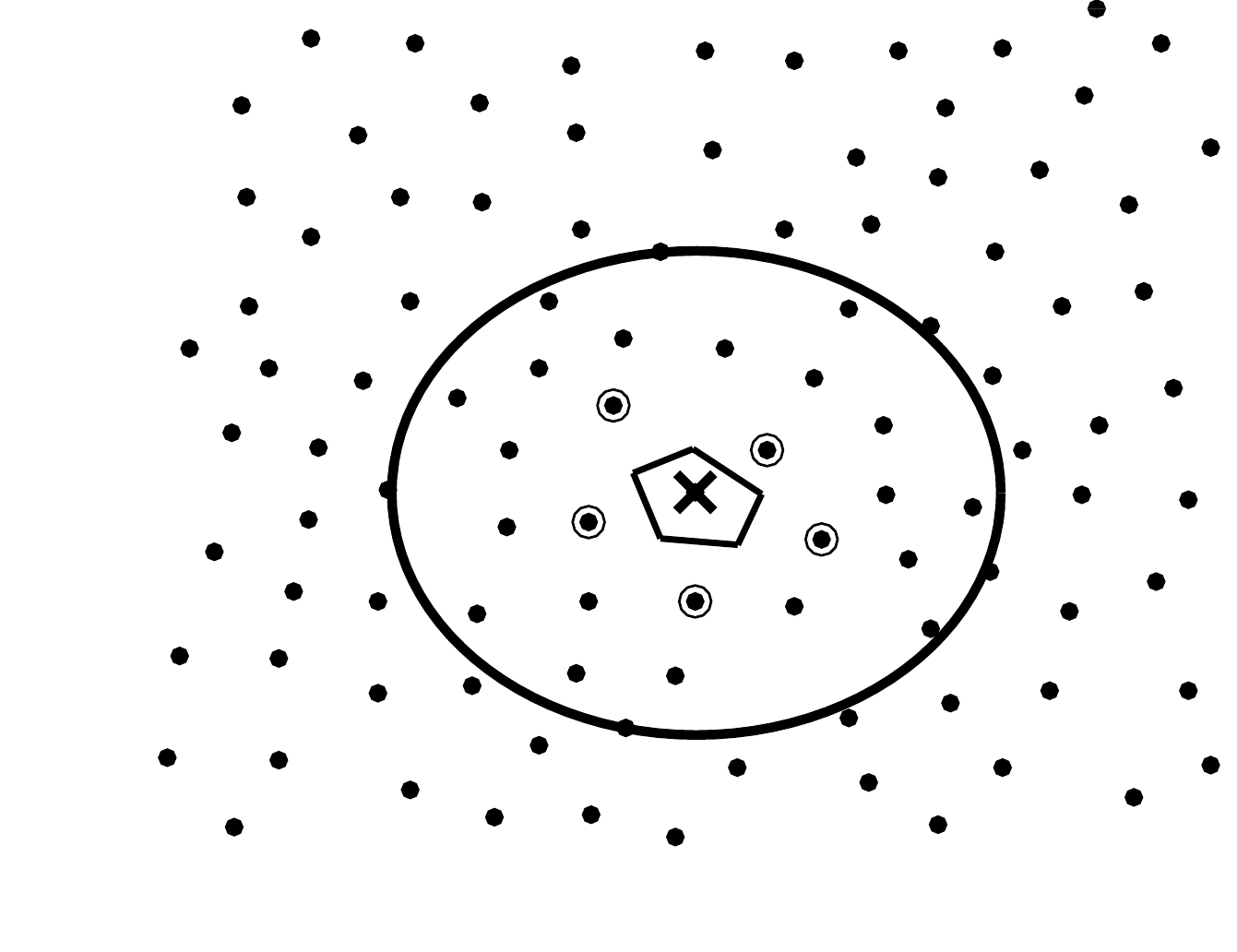}
  \includegraphics[width=0.4\textwidth]{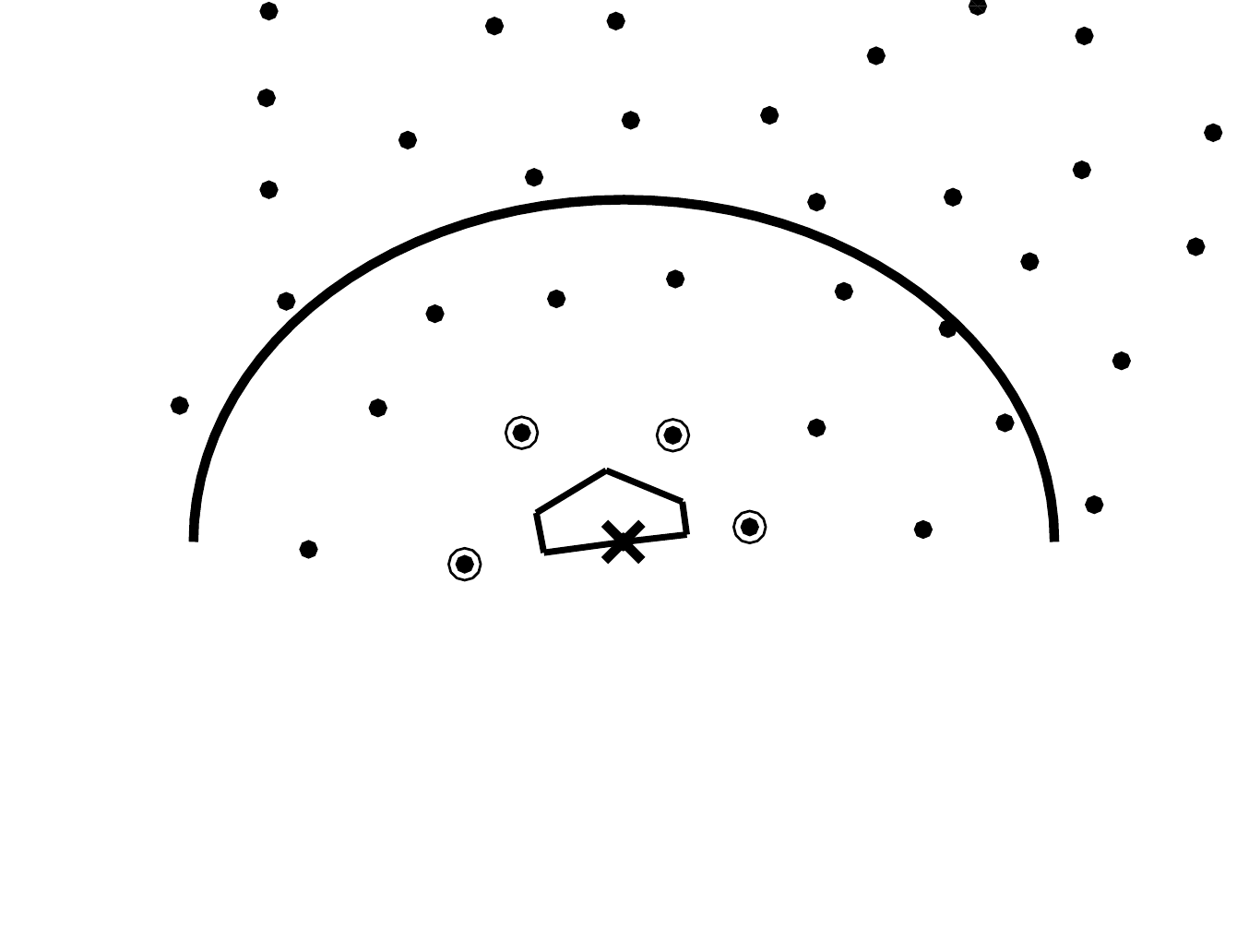}
  \caption{Control cell within a local support. Center point is marked with a cross. Nearest neighbours $l \in I_i$ are marked with circles. Interior point~(left) and boundary point~(right).}
  \label{Fig:ControlCellplusSupport}%
\end{figure}

For each point $i$, we consider the Delaunay tessellation of the $n_i$ points in its support domain $S_i$. The Voronoi diagram forms the dual graph to the Delaunay tessellation. Among the tessellations, the Voronoi cell containing point $i$ is the only one of interest, and is used as the control cell over which the flux balance is carried out.  Figure~\ref{Fig:ControlCellplusSupport} shows such a cell within the support domain. The point $i$ is associated with a volume $V_i$, taken to be the volume of this cell. $V_i$ is also used for accurate post-processing. The local tessellations can also be used in the computation of geometric parameters of the point cloud, such as the identification of free surface points in flows with open boundaries. Alternatively, instead of the Voronoi dual, a centroidal dual could also be used to the same effect.

Using the local tessellation for the support $S_i$, we obtain a set of nearest neighbours $I_i$, in the Delaunay sense. For $l \in  I_i$, we let $\widetilde{il}$ denote the dual edge~(in 2D) or face~(in 3D) of the Delaunay edge $il$. Further, we let $A_{il}$ denote the area of $\widetilde{il}$ and $\vec{n}_{il}$ denote the unit normal of $\widetilde{il}$, pointing away from $i$. For boundary points, the geometric area $A_i$ and outward pointing unit normal $\vec{n}_i$ are used to truncate and close the semi-infinite Voronoi cell, as shown in Figure~\ref{Fig:ControlCellplusSupport}. To simplify notation, we define the closure $\bar{I}_i$ of $I_i$ to include the point $i$ itself, if $i$ is a boundary point. Further, we let $\vec{n}_{ii} = \vec{n}_i$ and $A_{ii} = A_i$. 

\begin{figure}[!htbp]
  \centering
  \includegraphics[width=0.4\textwidth]{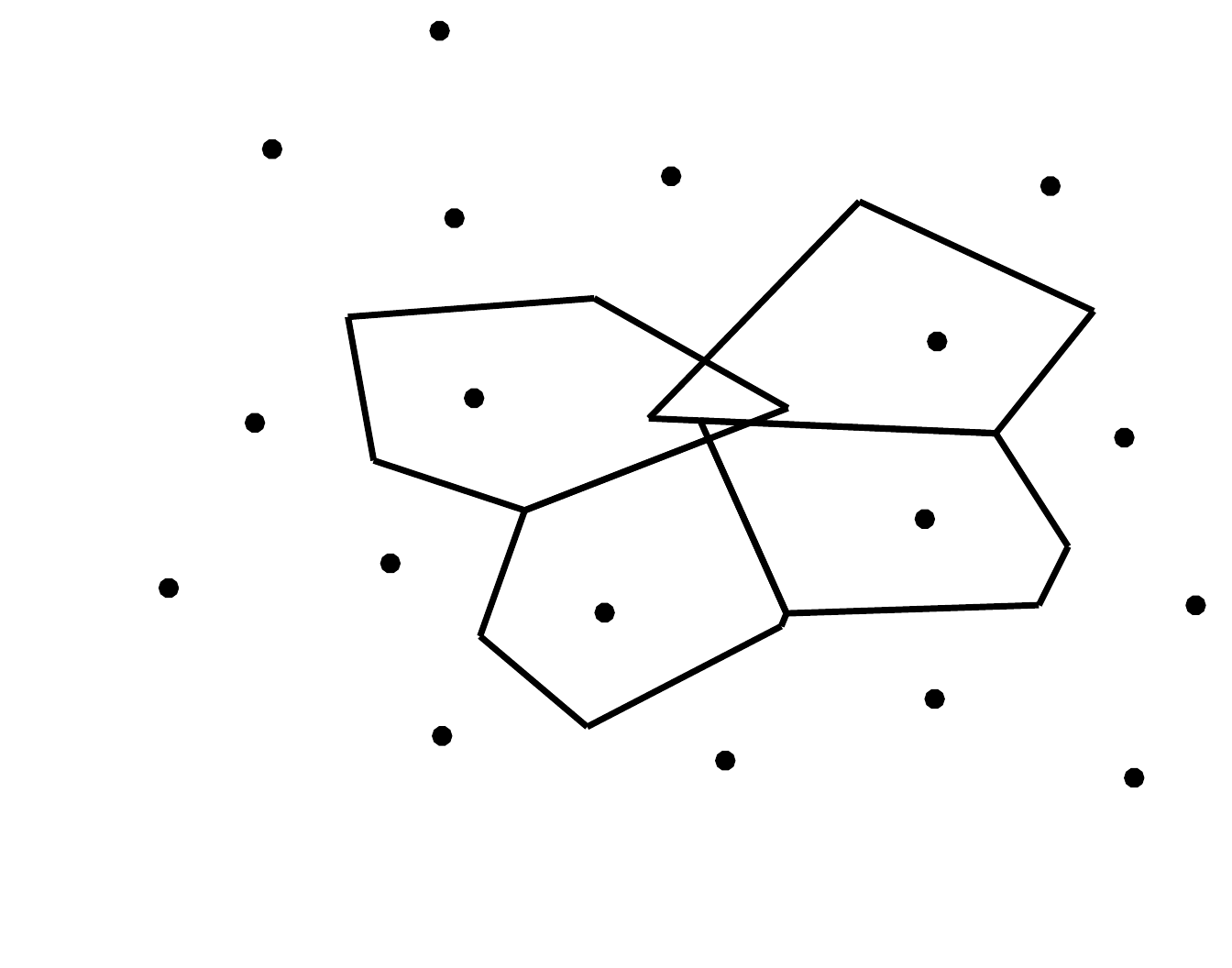}
  \phantom{abcdef}
  \includegraphics[width=0.4\textwidth]{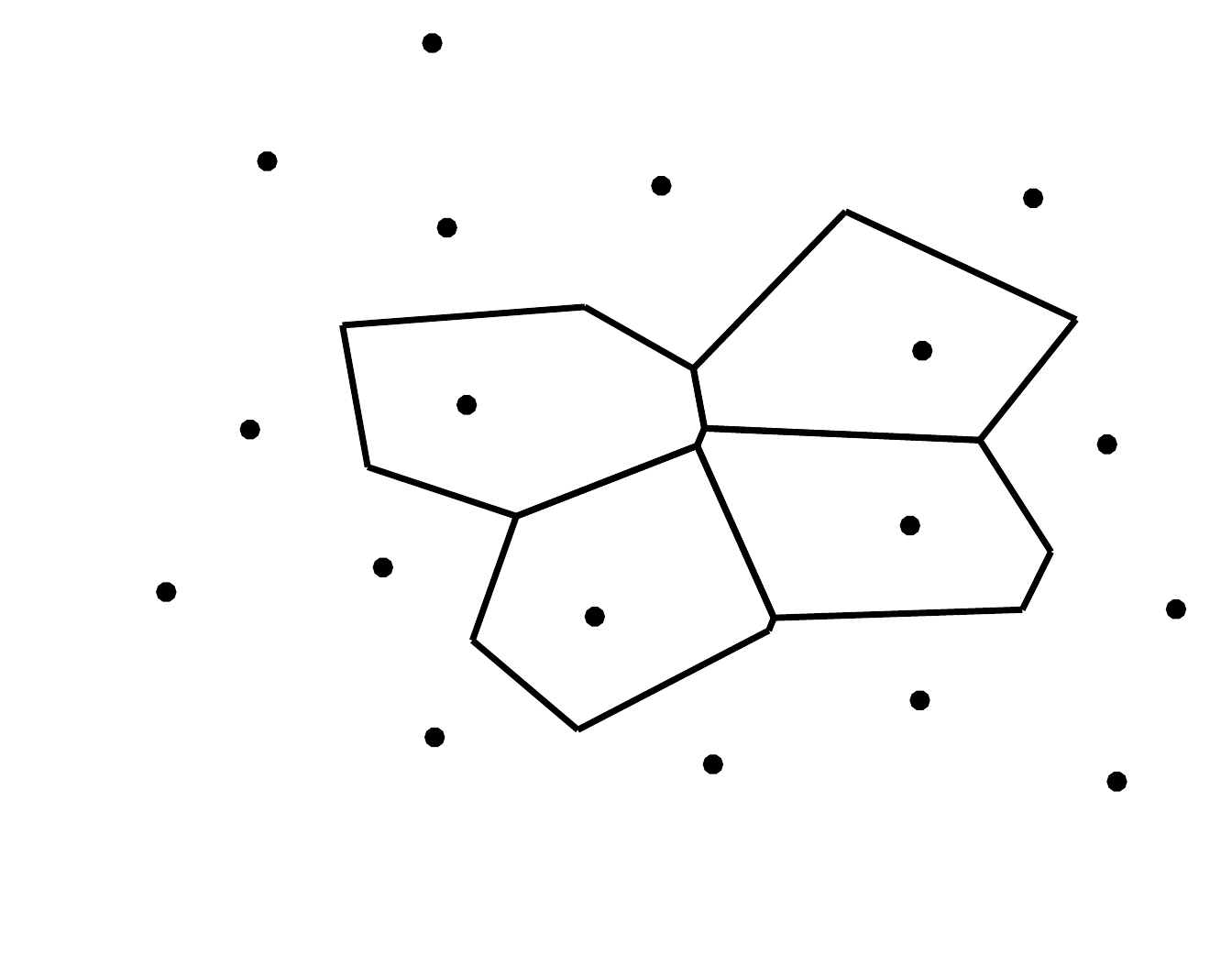}
  \caption{A 2D example of non-symmetric locally defined Voronoi cells~(left), and the same point configuration with symmetric locally defined Voronoi cells~(right). The non-symmetric cases occur when the support domains are not sufficiently large. The size of the support domains considered are $\beta = 0.6$, $mean(n_i)=11$~(left) and $\beta = 0.7$, $mean(n_i)=18$~(right).}
  \label{Fig:NonSymCells}%
\end{figure}

We note that since the Voronoi cells are defined locally on the support domain of each point, and not globally, they need not stitch together to form a global tiling of the computational domain. The uniqueness property of Delaunay tessellations suggests that the Voronoi cells should stitch together perfectly. However, that only holds for sufficiently large support domains. For small support domains, it is possible that an insufficient number of points in each neighbourhood leads to unsymmetric Voronoi cells. Under such situations, symmetry of the interfaces values, $A_{ij}$ and $A_{ji}$; $\vec n_{ij}$ and $\vec n_{ji}$, is violated. Figure~\ref{Fig:NonSymCells} shows such an example of adjacent non-symmetric cells. We measure the extent of stitching together to form a consistent global mesh by the error
\begin{equation}
	\label{Eq:MeshError}
	\epsilon_{mesh}= \frac{ \sum_{i=1}^N \sum_{l\in I_i} \|\vec{n}_{il}A_{il} + \vec{n}_{li}A_{li} \|  }{2\sum_{i=1}^N V_i}\,.
\end{equation}
If $i \notin I_l$, we set $\vec{n}_{li}=\vec{0}$ and $A_{li}=0$. If the cells are based on a global tessellation, $\vec{n}_{il}A_{il} = - \vec{n}_{li}A_{li}\,\forall i,l$, which leads to $\epsilon_{mesh}=0$. Table~\ref{tab:stitching} shows the extent of stitching together of locally defined Voronoi cells for different support sizes $\beta$, and the corresponding average number of points in the support domain. Table~\ref{tab:stitching} illustrates that for large enough support domains, which correspond to high values of $\beta$, the intersection of local Voronoi cells is minimum. Throughout this paper, we use support sizes in the range of $\beta\in[0.75,1]$. This results in approximately $45-55$ points in the support domain of each point in three dimensions, and about $20$ points in two dimensions.
\begin{table}[!htbp]
	\caption{Errors in formation of a global tiling by stitching together locally obtained Voronoi cells (in 3D). The domain considered is a unit sphere with $h=0.5$, $N=1400$.}
	\centering
	\tabsize
	\label{tab:stitching}
	\begin{tabular}{lll}
	\toprule
	$\beta$ & $\epsilon_{mesh}$ & $mean(n_i)$ \\
	\midrule
	$0.5$  &  $3.31$   & $16$ \\
	$0.55$ &  $1.08$   & $21$ \\	
	$0.6$  &  $0.33$   & $27$ \\
	$0.65$ &  $9\times 10^{-3}$   & $35$ \\
	$0.7$  &  $6\times 10^{-15}$  & $41$ \\
	\bottomrule
	\end{tabular}%
\end{table}

Several meshfree methods, such as EFG \cite{Belytschko1994} and PFEM \cite{Onate2004}, use a globally defined background mesh for a variety of purposes. Such a global background mesh could be used for the control cells, but defining cells based on local tessellations on a small number of points is significantly faster, especially when done in parallel. The concept of shapeless meshfree volumes for each point have also been proposed by several authors \cite{Chiu2012, Chiu2011,Katz2010,KatzThesis}, but their use presents several problems in the present context. Most importantly, the proposed meshless volumes are not closed, which, as we shall show later, is necessary for the construction of our flux conserving differential operators.
\subsection{Flux Conserving Differential Operators}
For the definition of our new numerical first derivative operators, consider the divergence theorem
\begin{equation}
	\label{Eq:DivTh_Vector}
	\int_{\Omega}\nabla\cdot\vec{g}\,dV = \int_{\partial\Omega}\vec{n}\cdot \vec{g} \,dA\,,
\end{equation}
for sufficiently smooth $\vec{g}$ and $\partial\Omega$. Setting $g$ to be $(f,0,0)$, $(0,f,0)$ and $(0,0,f)$ alternatingly, a scalar version of Eq.\,\eqref{Eq:DivTh_Vector} can be obtained as
\begin{equation}
	\label{Eq:DivTh2}
	\int_{\Omega}\partial^kf\,dV = \int_{\partial\Omega}fn^k\,dA\,.
\end{equation}
Throughout this paper, the superscript $k=x,y,z$ denotes the spatial dimension. A discrete local version of Eq.\,\eqref{Eq:DivTh2} on the control cells can be written as follows
\begin{equation}
	\label{Eq:LocalDDT}
	V_i\tilde{\partial}^k_if = \sum_{l\in \bar{I}_i} F_{il}(f)\,,
\end{equation}
where $F_{il}(f)$ is a numerical flux function which depends on $\vec{n}_{il}$, $A_{il}$, $f_i$ and $f_l$, and possibly their derivatives. One possibility for $F$, to achieve a local discrete divergence theorem, can be
\begin{equation}
	\label{Eq:FluxF}
	F_{il}(f) = 
	\begin{cases} 
		 f_in_i^kA_i & \mbox{if } i\in\partial\Omega \mbox{ and } l = i\,, \\
		 f_{il}n_{il}^k\,A_{il}    & \mbox{elsewhere} 		\,,
	\end{cases}
\end{equation}
with $f_{il}=\frac{f_i + f_l}{2}$. $\vec{n}_{il}$ and $A_{il}$ are determined locally based on the control cell as mentioned in Section\,\ref{sec:ControlCell}. Throughout this paper, we consider only simple symmetric flux functions similar to Eq.\,\eqref{Eq:FluxF}. However, more sophisticated flux functions can also be used and can also be coupled with flux functions in Eq.\,\eqref{Eq:FPM_DiffOp}. The use of flux functions in the definition of the differential operators, Eq.\,\eqref{Eq:FPM_DiffOp}, have been considered in \cite{Chiu2012, Chiu2011}; and those specific to the Finite Pointset Method have been studied in \cite{SeifarthThesis}.

For time-dependent problems, stencil coefficients for numerical differential operators are found such that they satisfy Eq.\,\eqref{Eq:LocalDDT} for conservation of fluxes at the previous time level, in addition to Eq.\,\eqref{Eq:Consistency}. Thus, if a field $\phi$ needs to be conserved, when proceeding from time-step $t^n$ to $t^{n+1}$, the discrete first-derivatives are found by the following minimization
\begin{align}
	\sum_{j\in S_i}c_{ij}^{k}m_j &= \partial^k_i m\qquad\ \forall m\in\mathcal{M}\,, \label{Eq:ConserStencil_consistency}\\
	\sum_{j\in S_i}c_{ij}^{k}\phi_j^{(n)} &= \frac{1}{V_i}\sum_{l\in \bar{I}_i}F_{il}(\phi^{(n)}) \,,  \label{Eq:ConserStencil}\\
	\text{min } J &= \sum_{j\in S_i} W_{ij}(c_{ij}^{k})^2\,, \label{Eq:ConserStencil_minimization}
\end{align}
where the bracketed superscript denotes the time-level and $k=x,y,z$. We emphasize that the superscript $k$ in Eq.\,\eqref{Eq:ConserStencil_consistency} - Eq.\,\eqref{Eq:ConserStencil_minimization} is used to denote only the first order derivatives in different directions, $x,y,z$, whereas the superscript ${}^*$ in Eq.\,\eqref{Eq:FPM_DiffOp} denotes any arbitrary differential operator. Eq.\,\eqref{Eq:ConserStencil_consistency} - Eq.\,\eqref{Eq:ConserStencil_minimization} constitute an enhancement of the GFDM differential operators by the conservation condition Eq.\,\eqref{Eq:ConserStencil}. If the monomials $m\in\mathcal{M}$ are taken up to second order, Eq.\,\eqref{Eq:ConserStencil_consistency} represents 10 conditions in 3D and 6 conditions in 2D. Whereas the number of unknown $c_{ij}$ values is the same as the number of points in each support domain, which is typically around 50 and 20, in 3D and 2D respectively, for the support sizes being used. Thus, there is plenty of numerical freedom to impose the additional flux conservation condition, Eq.\,\eqref{Eq:ConserStencil}, without affecting the accuracy of gradient reconstructions, Eq.\,\eqref{Eq:ConserStencil_consistency}. Clearly, the stencil coefficients $c_{ij}$ can be found locally at each point $i$, without the need for the global systems used by previous conservative meshfree GFDMs \cite{Chiu2012,Chiu2011}.

For each point $i=1,2,\dots,N$, Eq.\,\eqref{Eq:ConserStencil_consistency} represents the usual monomial consistency conditions as described by Eq.\,\eqref{Eq:Consistency}, which is dependent on all points in the support domain $S_i$. On the other hand, Eq.\,\eqref{Eq:ConserStencil} represents the conservation of numerical fluxes at the previous time-level, where the fluxes are defined solely on the locally defined control cell, which depends only on the nearest neighbours $\bar{I}_i$. This difference is illustrated in Figure~\ref{Fig:ControlCellplusSupport}. In explicit time-integration schemes where spatial derivatives are computed only at the previous time level, the differential operators defined by Eq.\,\eqref{Eq:ConserStencil_consistency} - Eq.\,\eqref{Eq:ConserStencil_minimization} degenerate to ones based solely on the locally defined control cell, making it resemble Voronoi-based FVMs. This is because the spatial derivatives used in such cases would be completely described by Eq.\,\eqref{Eq:ConserStencil}. Thus, the inclusion of Eq.\,\eqref{Eq:ConserStencil} is only done with implicit time-integration schemes, which are used throughout this paper.

The additional time taken due to the addition of Eq.\,\eqref{Eq:ConserStencil} is not significant. As mentioned earlier, tessellations are usually performed in classical GFDMs for prescribing volumes to points, areas to boundary points and the detection of free surfaces. Thus, no extra tessellations need to be performed to determine the local control cells. Secondly, the size of the systems being solved is changed by a very small amount. For example, in 3D, for second order-accurate differential operators, classical GFDMs solve a system of $10$ equations and about $50$ unknowns at each point. The addition of Eq.\,\eqref{Eq:ConserStencil} changes that to a system of $11$ equations and $50$ unknowns. Further, when the differential operators are defined locally as done in this paper, the largest portion of time of meshfree GFDM simulations is taken by the iterative solvers for the large sparse linear systems obtained by the discretization of the PDEs. The computation of the differential operators takes up a much smaller portion of the simulation time. Thus, the addition of the flux conservation constraint, Eq.\,\eqref{Eq:ConserStencil}, has a minimal effect on the overall time taken by the entire simulation. This is shown in the numerical examples in Section\,\ref{sec:3DChannel}.

Henceforth, for the sake of brevity, we denote this variant of GFDMs, in which the differential operators conserves fluxes of specific numerical fields, as FC-GFDM, where the FC stands for flux conserving. We note that Eq.\,\eqref{Eq:ConserStencil} does not ensure the conservation of numerical fluxes of the field determined at level $n+1$ with respect to the given differential operator stencils. Thus, even with respect to the field $\phi$, the use of differential operators defined by Eq.\,\eqref{Eq:ConserStencil_consistency}-Eq.\,\eqref{Eq:ConserStencil_minimization} would only form an approximately conservative method.

It should be ensured that the choice of the flux function $F$ is done such that Eq.\,\eqref{Eq:ConserStencil_consistency} and Eq.\,\eqref{Eq:ConserStencil} are consistent with each other when $\phi$ in the neighbourhood of $i$ is linearly dependent on the monomials $m\in\mathcal{M}$. For $F$ defined as in Eq.\,\eqref{Eq:FluxF}, this is ensured by the fact that the Voronoi cell is closed and not self-intersecting by definition. That guarantees that geometric fluxes are conserved. For $0$ order, the cell is closed, 
\begin{equation}
	\label{Eq:GeomCons0}
	\sum_{l \in \bar{I}_i} \vec{n}_{il} A_{il} = 0\,.
\end{equation}
Similarly, first order geometric conservation ensures
\begin{equation}
	\label{Eq:GeomCons1}
	\sum_{l \in \bar{I}_i}n_{il}^kx_{il}^{k'}A_{il} = \begin{cases} 0 &\mbox{if } k\neq k'\,, \\
	V_i & \mbox{if } k=k'\,, \end{cases}
\end{equation} 
where $k,k'= x,y,z$ denote the spatial dimension, $\vec{x}_{il}$ is the geometric center of the edge or face $\widetilde{il}$ for $l \neq i$, and $\vec{x}_{ii}=\vec{x}_i$. By the definition of the Voronoi diagram, $\vec{x}_{il}=\frac{\vec{x}_i + \vec{x}_l}{2}$. Eq.\,\eqref{Eq:GeomCons1} is also used to determine the volume of the cell.

\subsection{Higher Order Derivatives}
\label{sec:HigherOrderDerivatives}
The same procedure used above can be extended to compute numerical differential operators for higher order derivatives. This is illustrated with a diffusion operator. Consider the diffusion operator $D$ such that $D\phi = \nabla \cdot (\alpha \nabla \phi)$. The divergence theorem applied to the vector field $\alpha \nabla \phi$ leads to
\begin{equation}
	\label{Eq:Diff_DivTh}
	\int_\Omega \nabla \cdot (\alpha \nabla \phi)\,dV = \int_{\partial\Omega} \vec{n}\cdot (\alpha\nabla\phi)\,dA \,.
\end{equation}
A discrete local version of Eq.\,\eqref{Eq:Diff_DivTh} on the control cells can be written as
\begin{equation}
	\label{Eq:AdvDiff_LocConser_1}
	V_i\tilde{\partial}_i^D\phi = \sum_{l\in \bar{I}_i} G_{il}(\phi,\alpha) \,,
\end{equation}
where $G$ is a numerical flux function. Throughout this paper, $G$ is taken to be symmetric, similar to Eq.\,\eqref{Eq:FluxF}. 
\begin{equation}
	\label{Eq:FluxG}
	G_{il}(\phi) = 
	\begin{cases} 
		 \alpha_{i} \vec{n}_i \cdot\tilde{\partial}_{i}^{\nabla}\phi A_{i} & \mbox{if } i\in\partial\Omega \mbox{ and } l = i\,, \\
		 \alpha_{il} \vec{n}_{il} \cdot \tilde{\partial}_{il}^{\nabla}\phi\,A_{il}    & \mbox{elsewhere}\,,
	\end{cases}
\end{equation}
where $\alpha_{il}=\frac{\alpha_i + \alpha_l}{2}$, $\tilde{\partial}_{i}^{\nabla} = (\tilde{\partial}_{i}^{x};\tilde{\partial}_{i}^{y};\tilde{\partial}_{i}^{z})$ in 3D, and $\vec{n}_{il} \cdot \tilde{\partial}_{il}^{\nabla}$ is an approximation of the first order derivative in the direction of $\vec{n}_{il}$ and is given by
\begin{equation}
	\label{Eq:FluxG_2}
	\vec{n}_{il} \cdot \tilde{\partial}_{il}^{\nabla} \phi  = \frac{\phi_l-\phi_i}{\|\vec{x}_l-\vec{x}_i\|}\,.
\end{equation}
Once again, this method can easily be extended to use more sophisticated flux functions and less diffusive approximations than Eq.\,\eqref{Eq:FluxG_2}. Using Eq.\,\eqref{Eq:AdvDiff_LocConser_1}, when proceeding from time-step $t^n$ to $t^{n+1}$, the discrete diffusion operator is found by the following minimization
\begin{align}
	\sum_{j\in S_i}c_{ij}^{D}m_j &= \partial^D_i m\qquad\ \forall m\in\mathcal{M}\,, \label{Eq:NDiffusion1}\\
	\sum_{j\in S_i}c_{ij}^{D}\phi_j^{(n)} &= \frac{1}{V_i}\sum_{l\in \bar{I}_i}  G_{il}(\phi^{(n)} ,\alpha)\,, \label{Eq:NDiffusion2}\\
	\text{min } J &= \sum_{j\in S_i} W_{ij}(c_{ij}^{D})^2 \label{Eq:NDiffusion3}\,.
\end{align}
\vspace*{2\baselineskip}\\
The FC-GFDM differential operators can be used on both fixed and moving point clouds. However, their use on fixed point clouds has the significant disadvantage of needing differential operators to be recomputed at every time-step. On the other hand, since the traditional GFDM differential operators depend only on point locations, they can be computed in a single pre-processing step. This disadvantage is not present for moving point clouds, since changing point locations means that differential operators always need to be recomputed. We thus restrict the numerical study in the following sections to moving point clouds and Lagrangian frame of references.
\section{Advection Diffusion Equation}
\label{sec:AdvDiff}
Consider the advection-diffusion equation in Lagrangian formulation on a fixed domain $\Omega\subset\mathbb{R}^2$
\begin{align}
	\frac{D\vec{x}}{Dt} &= \vec{v}\,, \label{Eq:AdvDiffMovement}\\
	\frac{D\phi}{D t} &= \nabla \cdot (\alpha \nabla \phi) \label{Eq:AdvDiff}\,,
\end{align}
where $\phi$ is the concentration or temperature, $\alpha$ is the diffusivity, and $\vec{v}$ is the advection velocity which is assumed to be divergence-free. Energy conservation, Eq.\,\eqref{Eq:EnergyConservation}, leads to 
\begin{equation}
	\label{Eq:AdvDiff_EnergyCons}
	\frac{d}{d t}\int_{\Omega}\phi\,dV=\int_{\partial\Omega} \vec{n} \cdot \left(\alpha\nabla\phi - \vec{v}\phi \right) \,dA\,,
\end{equation}

\begin{figure}
  \centering
  \includegraphics[width=0.45\textwidth]{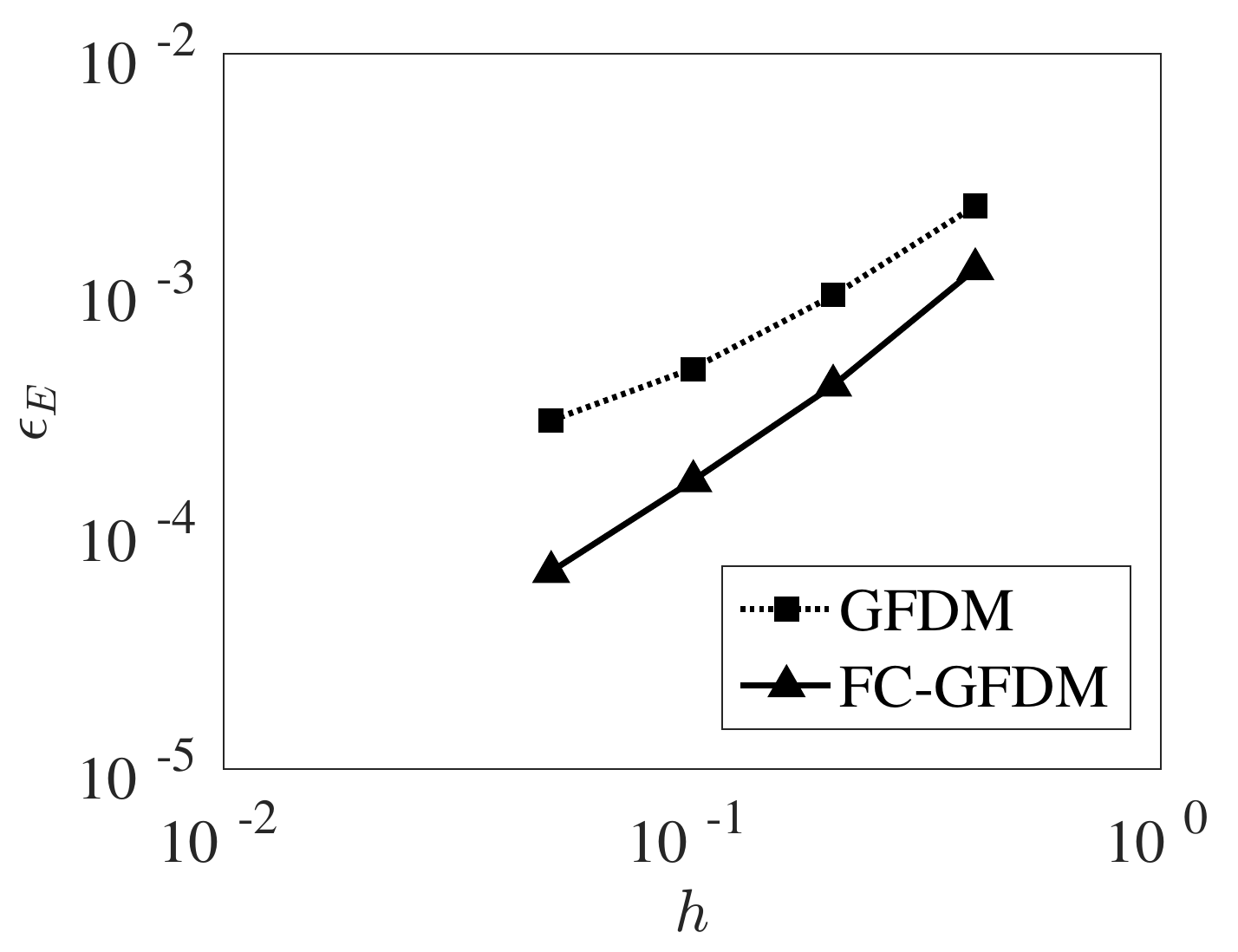}
  \caption{Convergence of error for the Advection Diffusion equation.}
  \label{Fig:EnergyChange_ADE}%
\end{figure}

For the FC-GFDM, the numerical diffusion operator are computed as done in Eq.\,\eqref{Eq:NDiffusion1}  - Eq.\,\eqref{Eq:NDiffusion3}. The spatial discretization of Eq.\,\eqref{Eq:AdvDiff} is obtained as
\begin{equation}
	\frac{D \phi_i}{D t} = \tilde{\partial}_i^{D}\phi = \sum_{j\in S_i}c_{ij}^{D}\phi_j \,.
\end{equation}	
%
In the first step, movement of the point cloud is performed in accordance with Eq.\,\eqref{Eq:AdvDiffMovement},
\begin{equation}
	\label{Eq:Movement}
	\vec{x}_i^{(n+1)} = \vec{x}_i^{(n)} + \vec{v}_i^{(n)}\Delta t + \frac{\vec{v}_i^{(n)} - \vec{v}_i^{(n-1)}}{\Delta t} (\Delta t)^2\,,
\end{equation}
for each point $i=1,2,\dots,N$. The explicit time-integration for the movement of points results in a CFL-like condition on the time-step size
\begin{equation}
	\label{Eq:VarDt}
	\Delta t= C_{\Delta t} \left( \frac{h}{\|\vec{v}\|} \right)_{min}\,.
\end{equation}
For the advection diffusion test case being considered, we use a coefficient of $C_{\Delta t} = 0.01$. The movement is followed by an implicit Euler time-integration of Eq.\,\eqref{Eq:AdvDiff}
\begin{equation}
	\label{Eq:AdvDiffDisc}
	\frac{\phi_i^{(n+1)} - \phi_i^{(n)}}{\Delta t} = \sum_{j\in S_i}c_{ij}^{D}\phi_j^{(n+1)} \,,
\end{equation}
The implicit system resulting from Eq.\,\eqref{Eq:AdvDiffDisc} can be solved using an iterative solver. In the results below, we use the BiCGSTAB solver \cite{BiCGSTAB}. Homogeneous Neumann boundary conditions are applied for $\phi$. Given these boundary conditions, error in energy conservation can be measured by
\begin{equation}
	\label{Eq:ADE_EnergyConsError}
	\epsilon_E = \frac{ \left| \int_{\Omega}\phi^{(end)}\,dV - \int_{\Omega}\phi^{(0)}\,dV  +  \int_0^{t_{end}}\left[ \int_{\partial\Omega} \vec{n} \cdot \vec{v}\phi \,dA \right]\,dt    \right| }{\int_{\Omega}\phi^{(0)}\,dV}\,,
\end{equation}
where the bracketed superscript $end$ denotes the values at the end of simulation, $t=t_{end}$, and the bracketed superscript $0$ denotes the initial condition. The computational domain $\Omega$ is taken to be $[-2,2]\times[-2,2]$. A uniform diffusion coefficient $\alpha=0.4$ and advection field $\vec{v}=(-y,x)$ are used. Initial conditions are taken to be
\begin{equation}
	\phi(\vec{x},0) = 
	\begin{cases} 
		 500 & \mbox{if } \|\vec{x}-(1,0)\|^2 < 0.1 \,, \\
		 0    & \mbox{elsewhere} 		\,.
	\end{cases} 
\end{equation}
The convergence of the error in energy conservation with a varying smoothing length is shown in Figure~\ref{Fig:EnergyChange_ADE}. The figure illustrates that the new flux conserving method produces much smaller errors. Both methods take similar simulation times. 
For an error $\epsilon$ and consecutive smoothing lengths $h_1$ and $h_2$, the rate of convergence of the solution with changing smoothing length is measured as 
\begin{equation}
	\label{Eq:ConvergenceRate}
	r = \frac{ \log \left( \frac{\epsilon(h_2)}{\epsilon(h_1)} \right) }{ %
	\log \left( \frac{h_2}{h_1} \right)}\,.
\end{equation}
The errors in energy conservation and the numerical orders of convergence of the errors are tabulated in Table~\ref{tab:AdvDiff}. 

\begin{table}
	\caption{Errors and experimental convergence rates for the advection diffusion test case. $h$ is the smoothing length, $N$ is the number of points in the entire domain at the initial state, $\epsilon_E$ is the error in energy conservation, and $r$ is the order of convergence of $\epsilon_{E}$.}
	\centering
	\tabsize
	\label{tab:AdvDiff}
	\begin{tabular}{l|c|c|c|c|c}
	\toprule
	&& \multicolumn{2}{c|}{GFDM} & \multicolumn{2}{|c}{FC-GFDM} \\
	$h$ & $N$ & $\epsilon_E$ & $r$ & $\epsilon_E$ & $r$ \\
	\midrule
	$0.4$  & $\phantom{11\,}752$  & $2.28\times 10^{-3}$ & $-$ & $1.24\times 10^{-3}$ & $-$\\
	$0.2$  & $\phantom{1}2\,778$ & $9.66\times 10^{-4}$ & $1.24$ & $4.05\times 10^{-4}$ & $1.61$\\
	$0.1$ & $10\,565$ & $4.76\times 10^{-4}$ & $1.02$ & $1.61\times 10^{-4}$ & $1.33$\\
	$0.05$& $36\,224$ & $2.89\times 10^{-4}$ & $0.72$ &  $6.73\times 10^{-5}$ & $1.26$\\
	\bottomrule
	\end{tabular}
\end{table}

\section{Incompressible Navier--Stokes Equations}
\label{sec:INSE}
Consider fluid flow modeled by the incompressible Navier--Stokes equations in Lagrangian form
\begin{align}
	\frac{D\vec{x}}{Dt} &= \vec{v} \,,\label{Eq:NSE1}\\
	\nabla\cdot\vec{v}  &= 0 \,,\label{Eq:NSE2}\\
	\frac{D\vec{v}}{D t} &= \frac{\eta}{\rho}\Delta \vec{v} - \frac{1}{\rho}\nabla p + \vec{g} \,,  \label{Eq:NSE3}
\end{align}
where $\vec{v}$ is the fluid velocity, $\rho$ is the density, $\eta$ is the dynamic viscosity and $\vec{g}$ includes both gravitational acceleration and body forces. We consider conservation with respect to the velocity field, which leads to
\begin{align}
	\int_\Omega \nabla\cdot \vec{v}\,dV &= \int_{\partial\Omega} \vec{n}\cdot \vec{v}\,dA \,,\\
	\int_\Omega \nabla\cdot (\vec{v}\otimes \vec{v})\,dV &= \int_{\partial\Omega} \vec{n}\cdot (\vec{v}\otimes \vec{v})\,dA \,, \\
	\int_\Omega \Delta\vec{v}\,dV &= \int_{\partial\Omega} \vec{n}\cdot\nabla \vec{v}\,dA \,.
\end{align}
For the FC-GFDM, the numerical Laplace operator is computed in a way similar to the numerical diffusion operator in Section~\ref{sec:HigherOrderDerivatives}.
\begin{align}
	\sum_{j\in S_i}c_{ij}^{\Delta}m_j &= \partial^{\Delta}_i m \qquad \forall m\in\mathcal{M} \,,\label{Eq:NSE_Lap1}\\
	\sum_{j\in S_i}c_{ij}^{\Delta}v_j^{k,(n)} &= \frac{1}{V_i}\sum_{l\in \bar{I}_i}  G_{il}(v_j^{k,(n)} ,1) \qquad k=x,y,z \,,\label{Eq:NSE_Lap2}\\
	\text{min } J &= \sum_{j\in S_i} W_{ij}(c_{ij}^{\Delta})^2 \label{Eq:NSE_Lap3}\,,
\end{align}
where the flux function $G$ is as defined earlier and $v_j^{k,(n)}$ denotes the $k$-component of the velocity at point $j$ and time-level $n$. The numerical first derivative approximations are computed as
\begin{align}
	\sum_{j\in S_i}c_{ij}^{k}m_j &= \partial^k_i m\qquad\ \forall m\in\mathcal{M} \,, \label{Eq:NSE_Grad1}\\
	\sum_{j\in S_i}c_{ij}^{k}v^{k,(n)}_j &= \frac{1}{V_i}\sum_{l\in \bar{I}_i}  F_{il}(v^{k,(n)}) \,, \label{Eq:NSE_Grad2}\\
	\sum_{j\in S_i}c_{ij}^{k}v^{k,(n)}_jv^{k',(n)}_j &= \frac{1}{V_i}\sum_{l\in \bar{I}_i}  F_{il}(v^{k,(n)}v^{k',(n)})\qquad k'=x,y,z \,, \label{Eq:NSE_Grad3}\\
	\text{min } J &= \sum_{j\in S_i} W_{ij}(c_{ij}^{k})^2 \,,\label{Eq:NSE_Grad4}
\end{align}
where the flux function $F$ is as defined earlier. While Eq.\,\eqref{Eq:NSE_Grad2} adds one extra constraint to the discrete differential operator definitions, Eq.\,\eqref{Eq:NSE_Grad3} adds two or three, depending on the spatial dimension.

Using the above defined numerical differential operators, Eq.\,\eqref{Eq:NSE1} - Eq.\,\eqref{Eq:NSE3} are solved with a meshfree projection method, similar to that of Chorin \cite{Chorin1968}. In the first step, the point cloud is moved by solving Eq.\,\eqref{Eq:NSE1} according to Eq.\,\eqref{Eq:Movement}. The projection method begins with the computation of an intermediate velocity $\vec{v}^*$ by solving the conservation of momentum equation, Eq.\,\eqref{Eq:NSE3},
\begin{equation}
	\label{Eq:vmm_vTild}
	\frac{\vec{v}^* - \vec{v}^{(n)}}{\Delta t} = \frac{\eta}{\rho}\Delta \vec{v}^* - \frac{1}{\rho}\nabla p^* + \vec{g}\,,
\end{equation}
where $p^*$ is a pressure guess found by updating the hydrostatic pressure, which is independent of the velocity \cite{Kuhnert2014}, and depends only on gravitational acceleration and body forces. $\vec{v}^*$ is then corrected by projecting it to a divergence-free space with the help of a correction pressure,
\begin{equation}
	\label{Eq:Projection}
	\vec{v}^{(n+1)} = \vec{v}^*-\frac{\Delta t}{\rho}\nabla p_{corr}\,.
\end{equation}
$p_{corr}$ is found by applying the divergence operator to Eq.\,\eqref{Eq:Projection} to obtain the pressure-Poisson equation
\begin{equation}
	\label{Eq:PressurePoisson}
	\frac{\Delta t}{\rho}\Delta p_{corr} = \nabla\cdot\vec{v}^*\,.
\end{equation}
The final pressure is then updated
\begin{equation}
	p^{(n+1)} = p^* + p_{corr}\,.
\end{equation}
Such meshfree projection methods have been widely used. For further details, see, for example, \cite{Jefferies2015,Kuhnert2014,Tiwari2016,Tiwari2002}. The spatial discretization of Eq.\,\eqref{Eq:vmm_vTild} - Eq.\,\eqref{Eq:PressurePoisson} are done using the differential operators defined earlier in this section, Eq.\,\eqref{Eq:NSE_Lap1} - Eq.\,\eqref{Eq:NSE_Lap3}, and Eq.\,\eqref{Eq:NSE_Grad1} - Eq.\,\eqref{Eq:NSE_Grad4}. The implicit systems obtained from Eq.\,\eqref{Eq:vmm_vTild} and Eq.\,\eqref{Eq:PressurePoisson} are solved using the BiCGSTAB iterative solver \cite{BiCGSTAB} without the use of any preconditioner. 
\subsection{Decaying Shear Flow}
\label{sec:DecayingShear}
\begin{figure}
  \centering
  \includegraphics[width=0.45\textwidth]{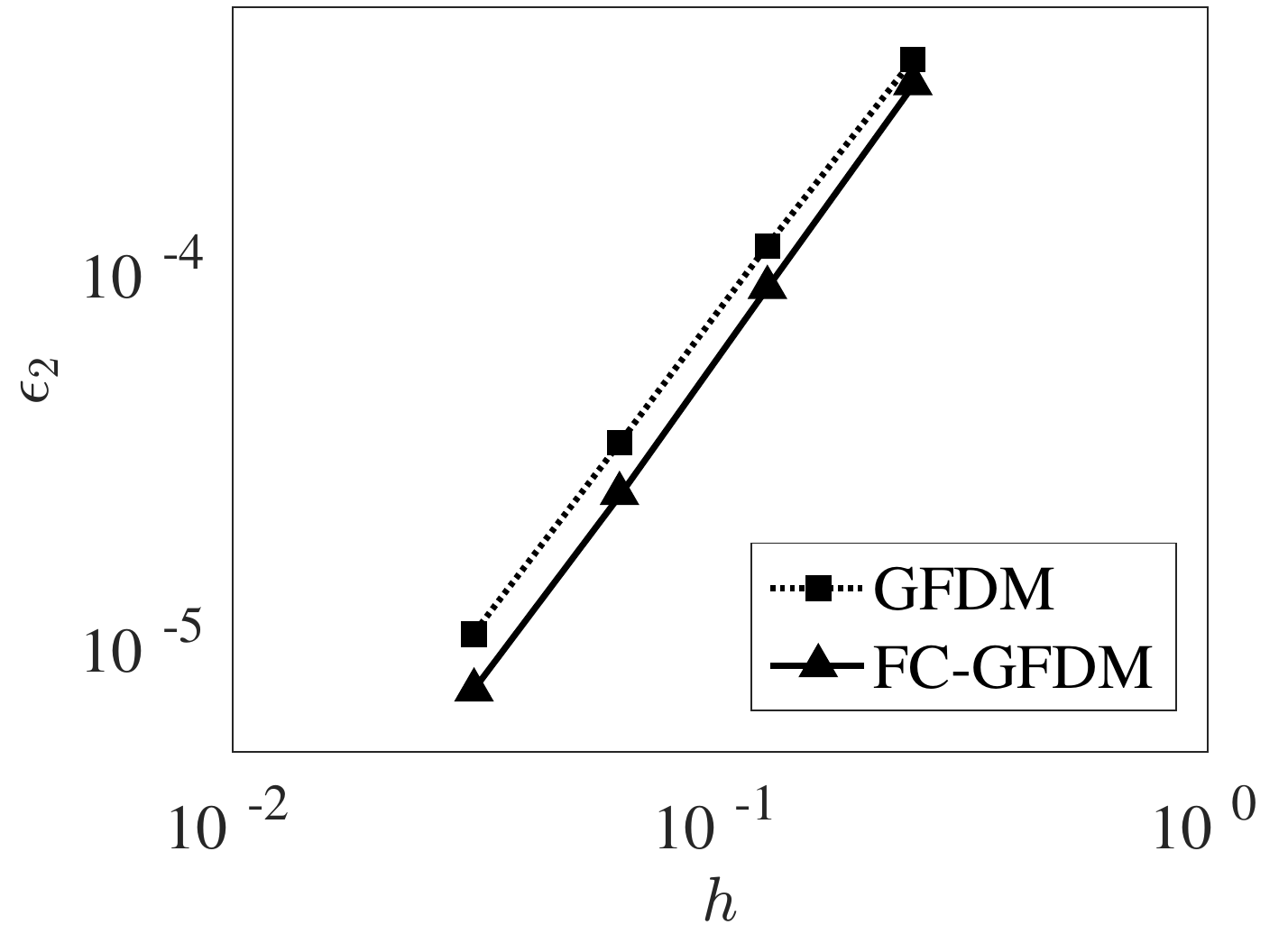}
  \caption{Errors for the decaying shear flow test case.}%
  \label{Fig:DecayingShear}
\end{figure}

For a validation case, we consider the case of decaying shear flow as reported in \cite{Guo2000}. For unit density, $\rho = 1 kg/m^3$, 
the analytical solution is given by
\begin{align}
	u_{\text{exact}} &= 1 \,,\\
	v_{\text{exact}} &= \cos(x-t)\exp(-\eta t) \,,\\
	p_{\text{exact}} & = p(t)\,,\\
	\vec{g} & = 0\,,
\end{align}
where $\vec{v}_{\text{exact}} = (u_{\text{exact}},v_{\text{exact}})$ is the analytical solution for the velocity, and $p_{\text{exact}}$ is the analytical solution for the pressure. The pressure is taken as $p(t) = 3.0 + 0.01\sin(20\pi t)$, as done in \cite{Guo2000}. Initial and boundary conditions are set in accordance with the analytical solution. The computational domain is taken to be $[0,1]\times[0,1]$. The time step is determined according to Eq.\,\eqref{Eq:VarDt} with a small coefficient of $C_{\Delta t} = 0.005$. The relative errors of the numerical solution are measured as 
\begin{equation}
	\label{Eq:PF_Error_2}
	\epsilon_2 = \left[ \frac{\sum \|\vec{v}_{\text{num}} - \vec{v}_{\text{exact}}\|^2}{\sum \|\vec{v}_{\text{exact}}\|^2}	\right]^{\frac{1}{2}} \,, 
\end{equation}
where $\vec{v}_{\text{num}}$ is the numerically obtained velocity. The errors in the solutions at $t = 1s$  with and without the addition of the flux conservation condition are shown in Figure~\ref{Fig:DecayingShear}. The figure illustrates that both methods show very similar results, with FC-GFDM being slightly more accurate. The errors are also tabulated, along with numerical convergence orders according to  Eq.\,\eqref{Eq:ConvergenceRate}, in Table~\ref{tab:DecayingShear}. Experimental orders of convergence are slightly higher for the new FC-GFDM. 

\begin{table}
	\caption{Errors and experimental convergence rates for the decaying shear flow test case. $h$ is the smoothing length, $N$ is the number of points in the entire domain at the initial state, $\epsilon_2$ is the error in energy conservation, and $r$ is the order of convergence of $\epsilon_{2}$.}
	\centering
	\tabsize
	\label{tab:DecayingShear}
	\begin{tabular}{l|c|c|c|c|c}
	\toprule
	&& \multicolumn{2}{c|}{GFDM} & \multicolumn{2}{|c}{FC-GFDM} \\
	$h$ & $N$ & $\epsilon_2$ & $r$ & $\epsilon_2$ & $r$ \\
	\midrule
	$0.25$  & $\phantom{1\,}161$  & $3.65\times 10^{-4}$ & $-$ & $3.11\times 10^{-4}$ & $-$\\
	$0.125$  & $\phantom{1\,}493$ & $1.15\times 10^{-4}$ & $1.67$ & $8.80\times 10^{-5}$ & $1.82$\\
	$0.0625$ & $1\,704$ & $3.41\times 10^{-5}$ & $1.75$ & $2.46\times 10^{-5}$ & $1.84$\\
	$0.03125$ & $6\,293$ & $1.04\times 10^{-5}$ & $1.71$ &  $7.36\times 10^{-6}$ & $1.74$\\
	\bottomrule
	\end{tabular}
\end{table}

For ``simple" domains and Dirichlet boundary conditions on all boundaries, the FC-GFDM shows only minor improvements over the classical GFDM. However, as the results in the coming sections show, FC-GFDM makes significant improvements for more complex simulation domains and boundary conditions.
\subsection{Flow Past a Square Cylinder}
\label{sec:2DChannel}
\begin{figure}[!htbp]
  \centering
  \includegraphics[width=0.45\textwidth]{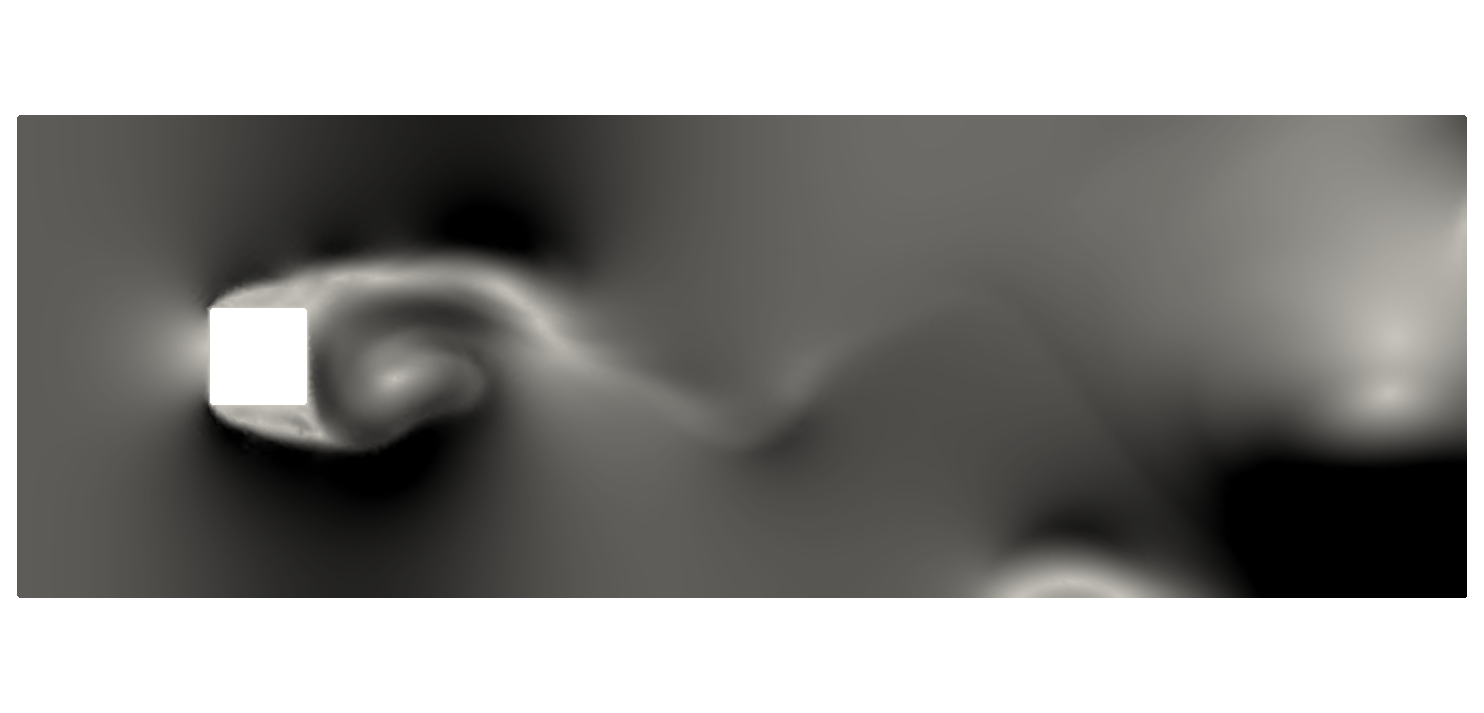}
  \includegraphics[width=0.45\textwidth]{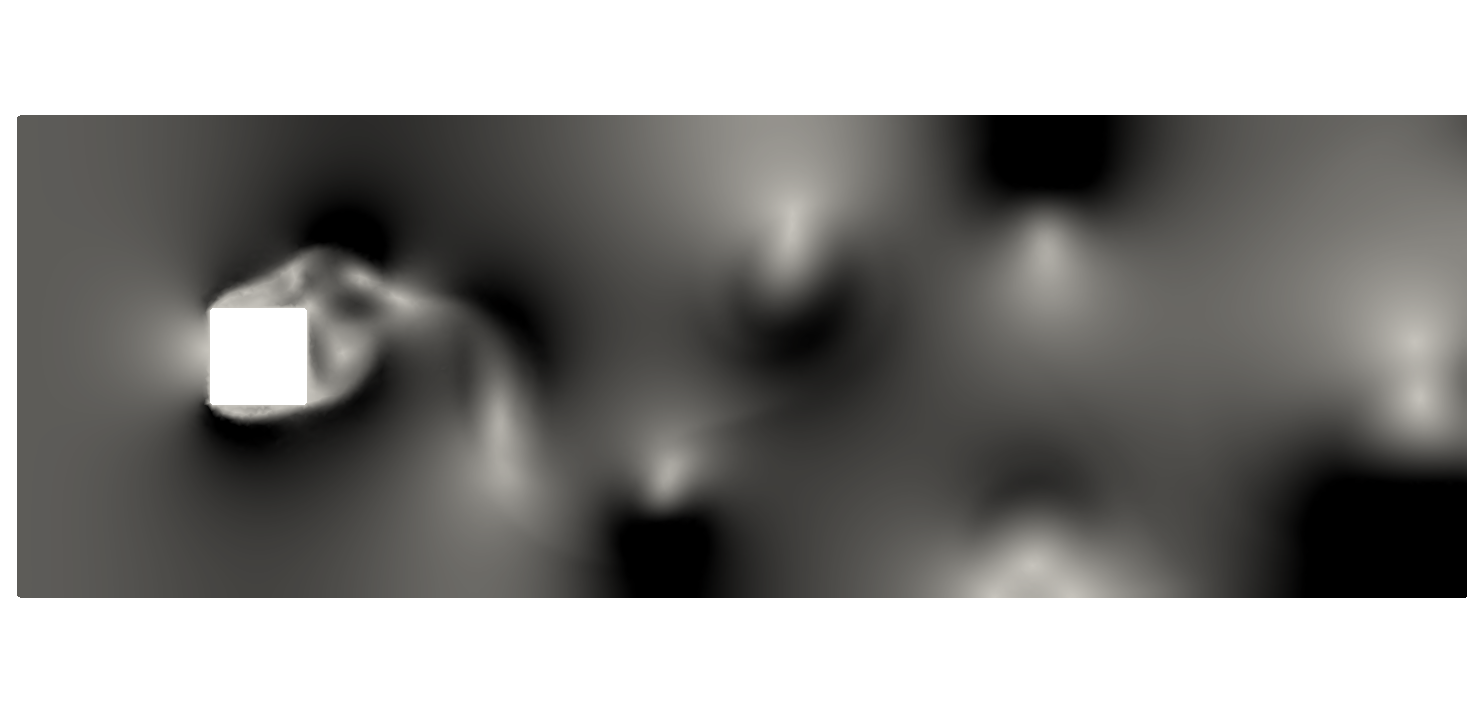}
  \caption{Flow past a square cylinder. Velocity profile in an inverted color scale for $Re = 10000$. $t=23s$~(left) and $t=34s$~(right).}%
  \label{Fig:2DChannel}
\end{figure}

We consider the flow past a square cylinder in two spatial dimensions as illustrated in Figure~\ref{Fig:2DChannel}, which has been a widely studied problem (see, for example, \cite{Saha2003, Sen2011}). 

The parameters of the simulations are $\rho = 1 kg/m^3$ and $t_{end} = 50s$. The viscosity is given by
\begin{equation}
	\eta = \frac{\|\vec{v}_{in}\| L}{Re} \,,
\end{equation}
where $\vec{v}_{in}=(2,0)m/s$ is the constant inflow velocity which is also used as the initial condition, $L = 30m$ is the length of the channel and $Re$ is the Reynolds number of the flow. Homogeneous Neumann boundary conditions for the velocity are used at the outflow and no-slip conditions on the walls. The pressure is kept constant at atmospheric pressure at the outflow and homogeneous Neumann boundary conditions are considered elsewhere for the pressure. A varying time-step is used according to Eq.\,\eqref{Eq:VarDt} with $C_{\Delta t}=0.3$, which results in larger time-steps than those considered in the earlier examples. The error in mass conservation is measured as the difference between the total volume of fluid flowing in and that flowing out, throughout the entire simulation, and is given by
\begin{equation}
	\epsilon_{mass} = \left| \frac{\int_0^{t_{end}} \left[\int_{ \partial \Omega_{in}}\vec{n}\cdot\vec{v}\,dA\right] dt + \int_0^{t_{end}} \left[ \int_{\partial\Omega_{out}}\vec{n}\cdot\vec{v}\,dA\right] dt }{\int_0^{t_{end}} \left[ \int_{\partial\Omega_{in}}\vec{n}\cdot\vec{v}\,dA\right] dt} \right| \,.
\end{equation}
The error in a global discrete divergence theorem at the new time level is given by
\begin{equation}
	\label{Eq:ErrorLDDT}
	\epsilon_{DDT}	= \frac{\left| \int_\Omega \nabla\cdot \vec{v}\, dV - \int_{\partial\Omega} \vec{n}\cdot \vec{v} \, dA \right| } {\int_{\Omega} \, dV } \,.
\end{equation}

\begin{figure}[!htbp]
  \centering
  \includegraphics[width=0.45\textwidth]{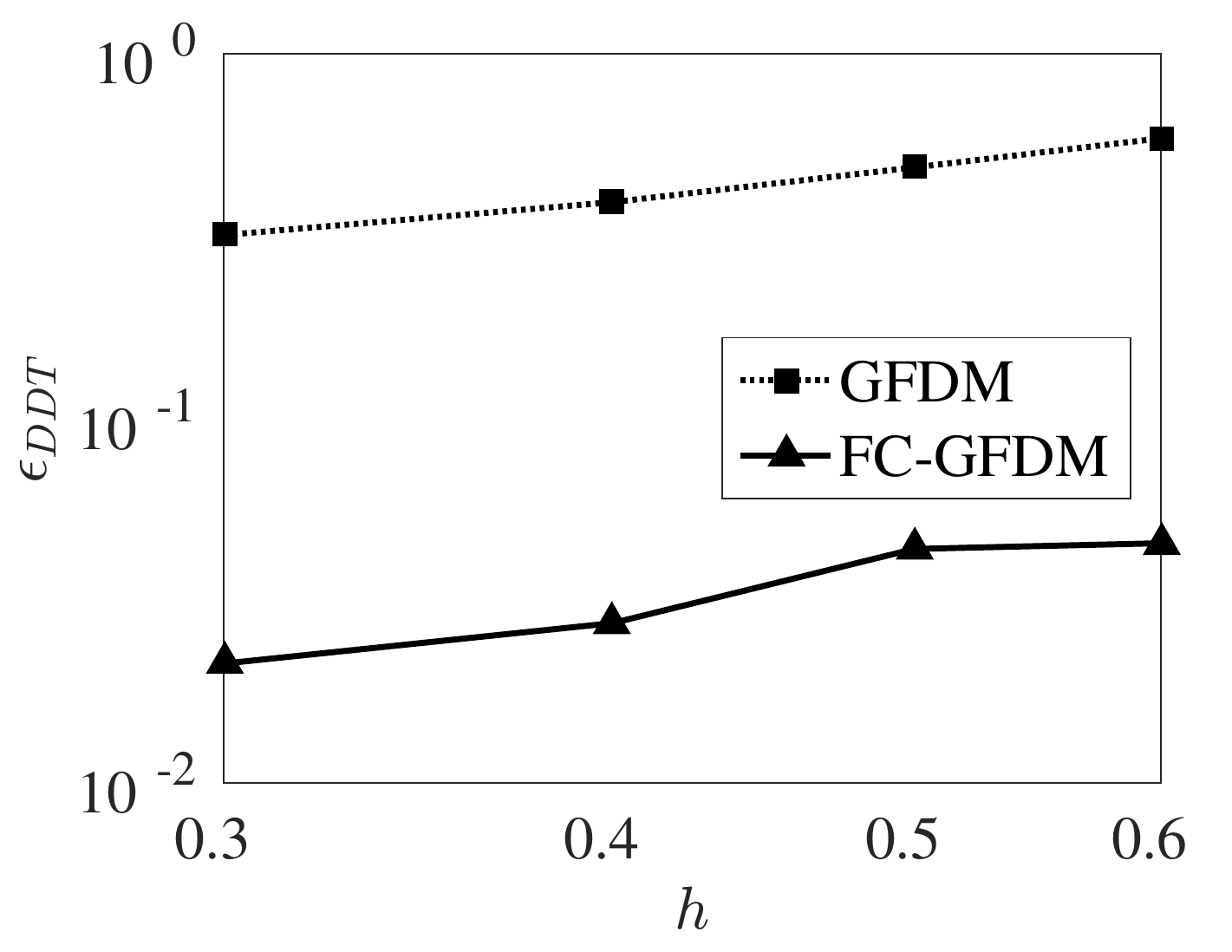}
  \includegraphics[width=0.45\textwidth]{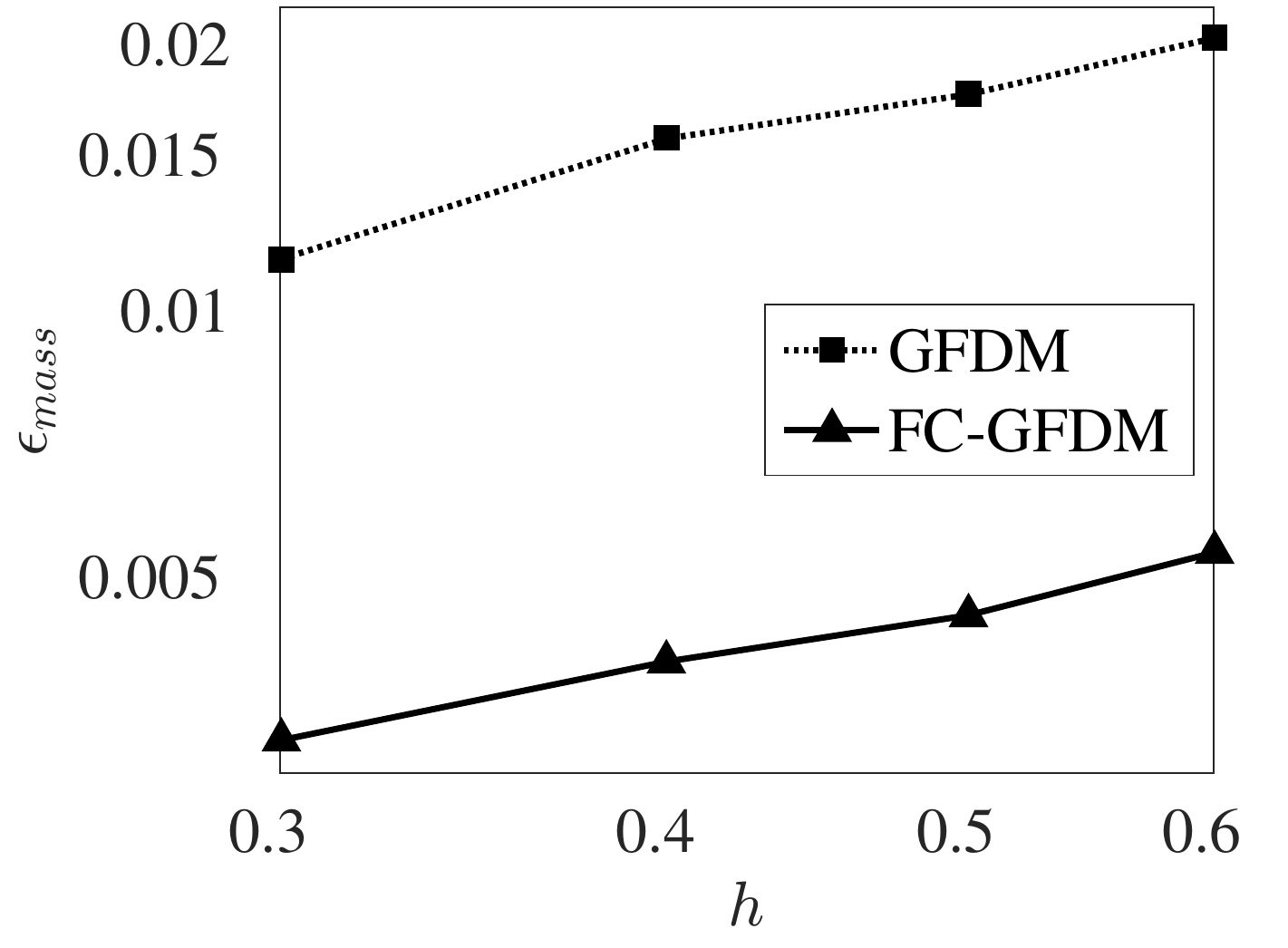}
   \caption{Flow past square cylinder: convergence of errors with smoothing length $h$ on log-log plots. Error in discrete divergence theorem~(left) and error in mass conservation~(right). The total number of points in the computational domain for the different smoothing lengths considered are shown in Table~\ref{tab:hvsN_2DChannel}.}
  \label{Fig:2DChannel_hConvergence}
\end{figure}

For a Reynolds number of $Re=10^4$, the convergence with respect to the smoothing length $h$ of $\epsilon_{DDT}$ averaged over all time steps and $\epsilon_{mass}$ for the two methods are shown in Figure~\ref{Fig:2DChannel_hConvergence}. The errors in the global discrete divergence theorem for the new FC-GFDM are lower by one order of magnitude when compared to that in the case of GFDM. Correspondingly, the errors in mass conservation are also lower by one order of magnitude in the case of FC-GFDM. The total number of points in the computational domain at $t=0s$ for the different smoothing lengths considered are shown in Table~\ref{tab:hvsN_2DChannel}.
\begin{table}
	\caption{$h$ vs. $N$ for flow past square cylinder}
	\centering
	\tabsize
	\label{tab:hvsN_2DChannel}
	\begin{tabular}{ll}
	\toprule
	Smoothing Length $h$ & Total number of initial Points $N$ \\
	\midrule
	$0.6$  &  $\phantom{1}3\,611$  \\
	$0.5$  &  $\phantom{1}5\,090$  \\
	$0.4$  &  $\phantom{1}7\,822$  \\
	$0.3$  &  $13\,661$  \\
	\bottomrule
	\end{tabular}
\end{table}

\begin{figure}
  \centering
  \includegraphics[width=0.45\textwidth]{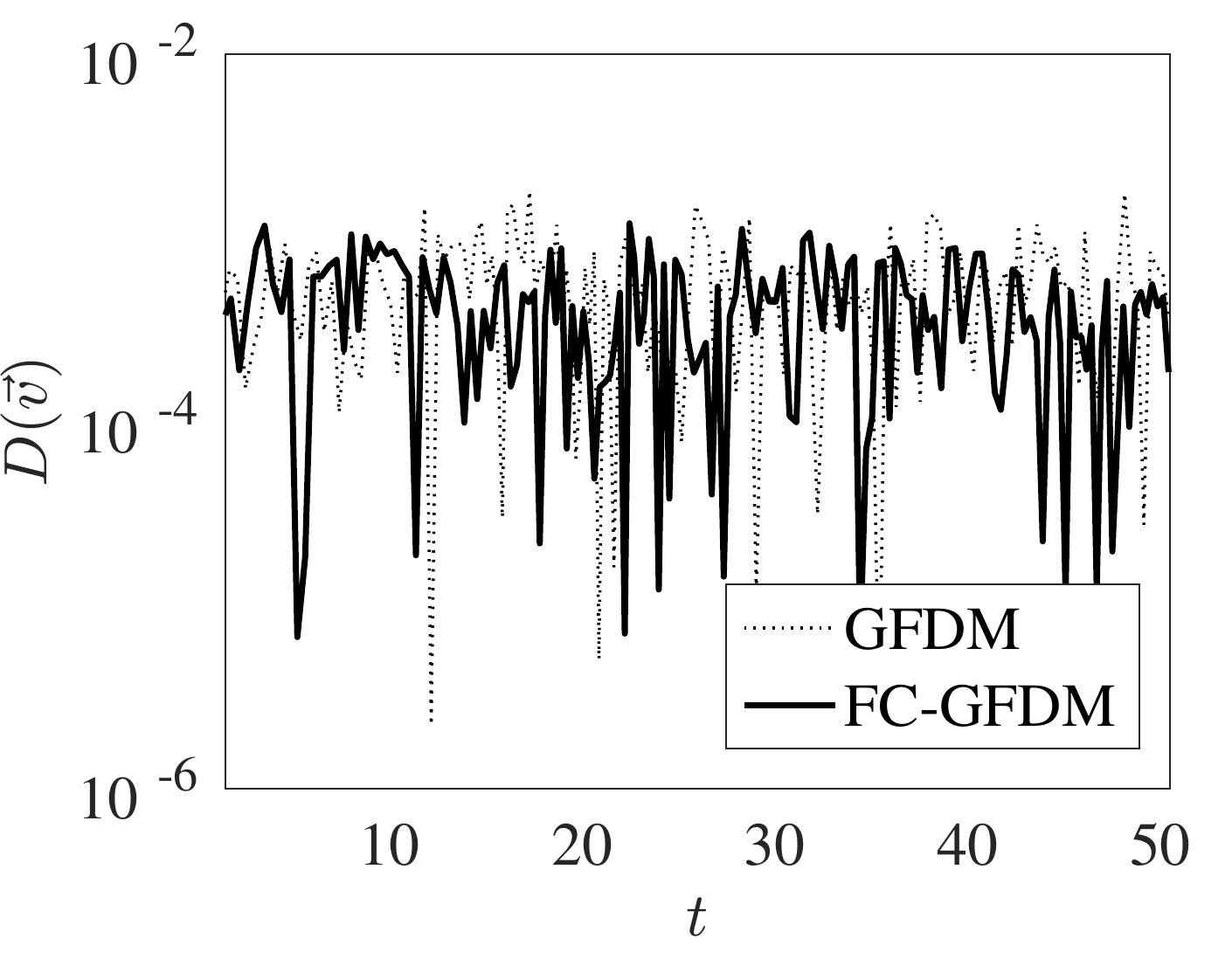}
  \includegraphics[width=0.45\textwidth]{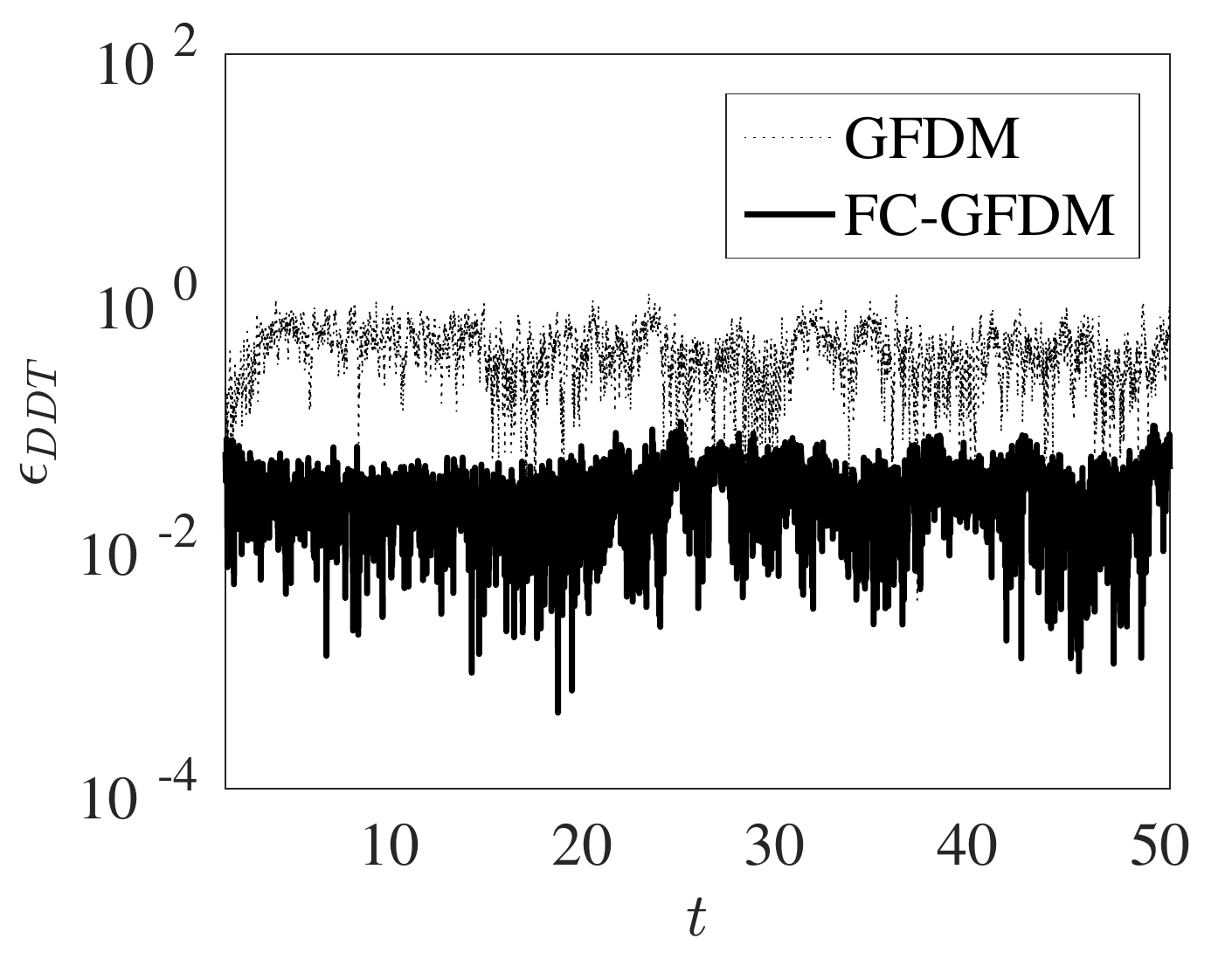}
  \includegraphics[width=0.45\textwidth]{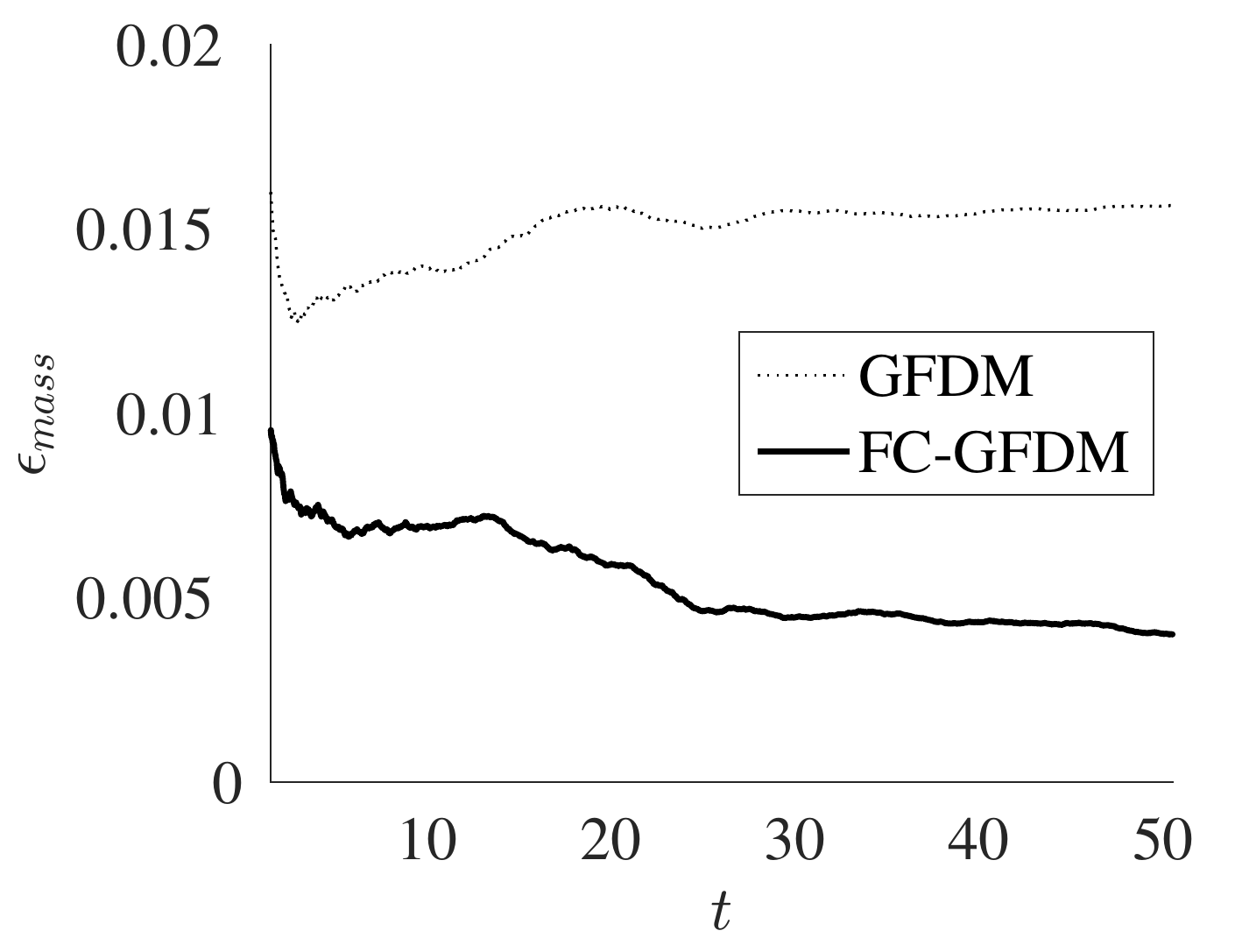}
   \caption{Flow past square cylinder $Re=10000$, $h=0.4$: Time evolution of total divergence~(top left), error in discrete divergence theorem~(top right) and error in mass conservation~(bottom). 1 in every 40 time steps is plotted in the divergence figure, whereas all time steps are plotted in the other two figures.}
  \label{Fig:2DChannel_TimeEvolution}
\end{figure}

A similar trend is also observed for lower Reynolds numbers of $10^3$ and $500$, but the errors are smaller for smaller Reynolds flow. The time evolution of the errors for $h= 0.4$ is shown in Figure~\ref{Fig:2DChannel_TimeEvolution}. This also includes a measure for the total divergence in the domain
\begin{equation}
	\label{Eq:Divergence}
	D(\vec{v}) = \frac{\int_{\Omega} |\nabla \cdot\vec{v}|\, dV}{\int_{\Omega} dV}\,.
\end{equation}
Figure~\ref{Fig:2DChannel_TimeEvolution} illustrates that the FC-GFDM does not significantly affect the velocity divergence values, and thus, the improvement in mass conservation is solely due to smaller errors in a global discrete divergence theorem.

\subsection{3D Obstructed Channel Flow}
\label{sec:3DChannel}
\begin{figure}
	\centering
	\includegraphics[scale=0.4]{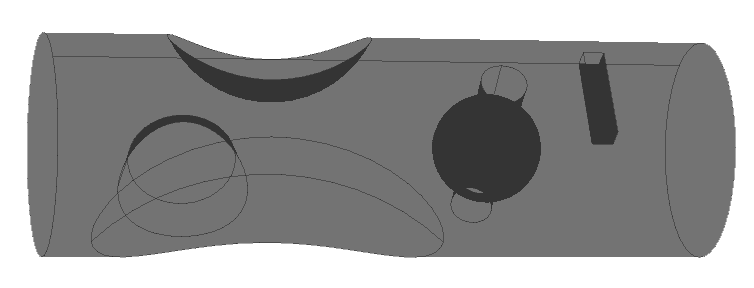}
	\caption{3D Channel with multiple obstructions. Fluid inflow is on the left, and outflow on the right}
	\label{Fig:3DChannel}%
\end{figure}

The case considered in the previous section is extended to three dimensions, with multiple obstructions in the channel, both convex and concave. The domain considered is shown in Figure~\ref{Fig:3DChannel}. 
For the simulations, we consider $\rho = 10^3 kg/m^3$, $\eta = 10^5 Pa\,s$, and $t_{end} = 1s$. The velocity at the inflow, on the left of the tube, is kept constant at $\vec{v}_{in}=(0,0,2)m/s$. The remaining boundary conditions are set up exactly as done in the 2D case mentioned earlier. The length of the channel is $6m$ which results in a Reynolds number in the order of magnitude of $10^{-1}$. Measurement of errors are also done as in the previous section.

\begin{figure}
  \centering
  \includegraphics[width=0.45\textwidth]{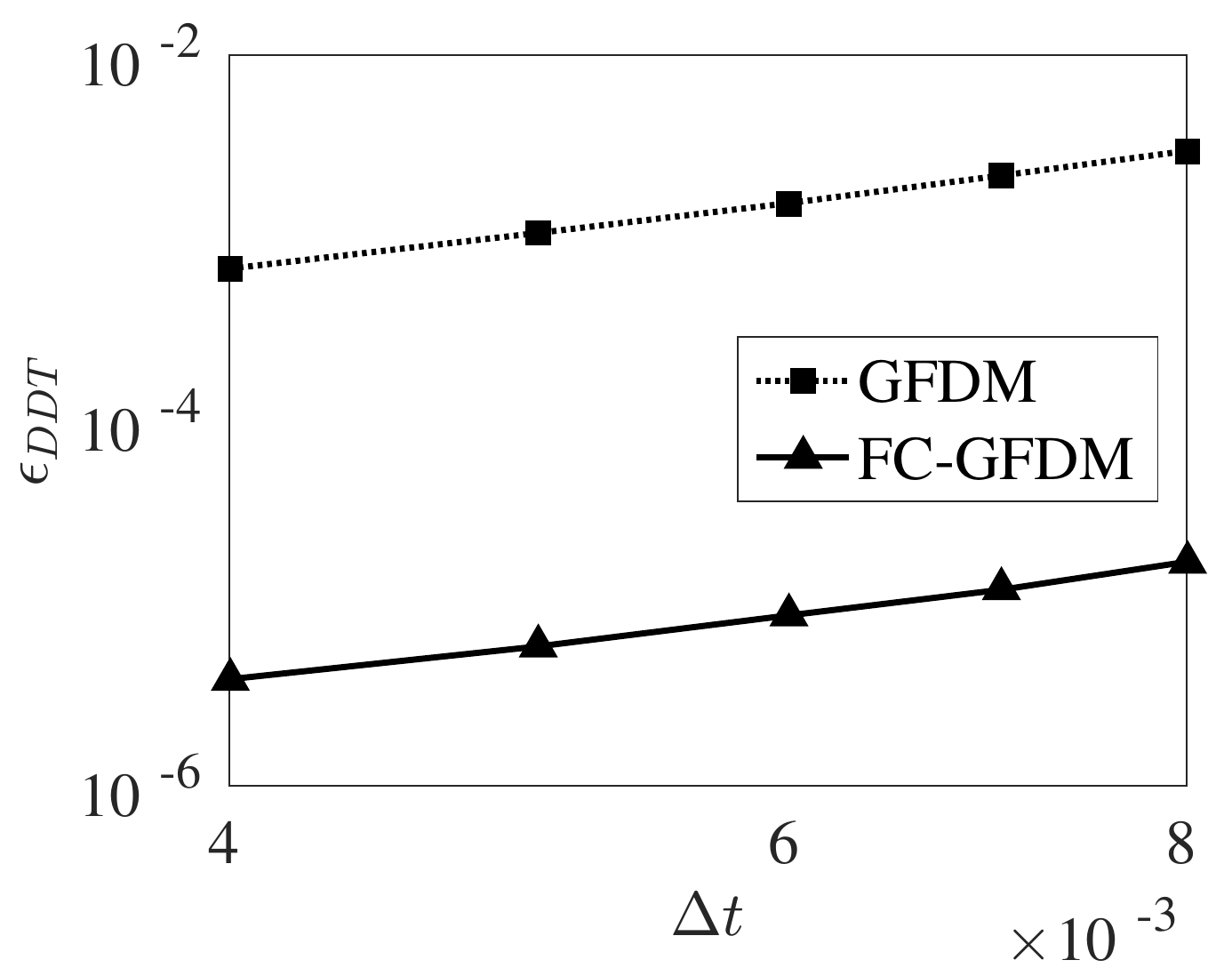}
	\includegraphics[width=0.45\textwidth]{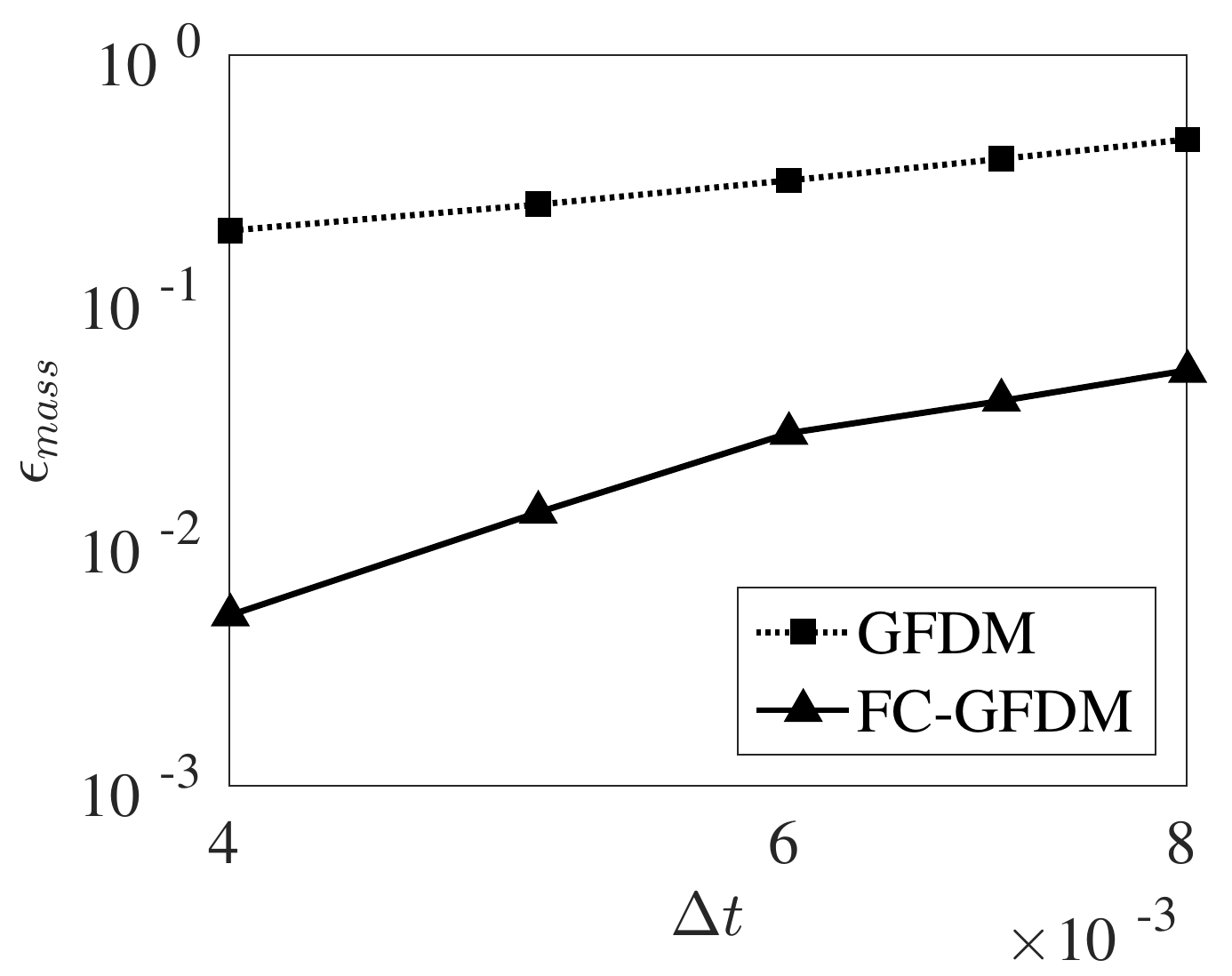}
  \caption{3D obstructed channel: convergence of solution with time step $\Delta t$. Error in discrete divergence theorem~(left) and error in mass conservation~(right)}
  \label{Fig:3DChannel_dt}%
\end{figure}

\begin{figure}
  \centering
  \includegraphics[width=0.45\textwidth]{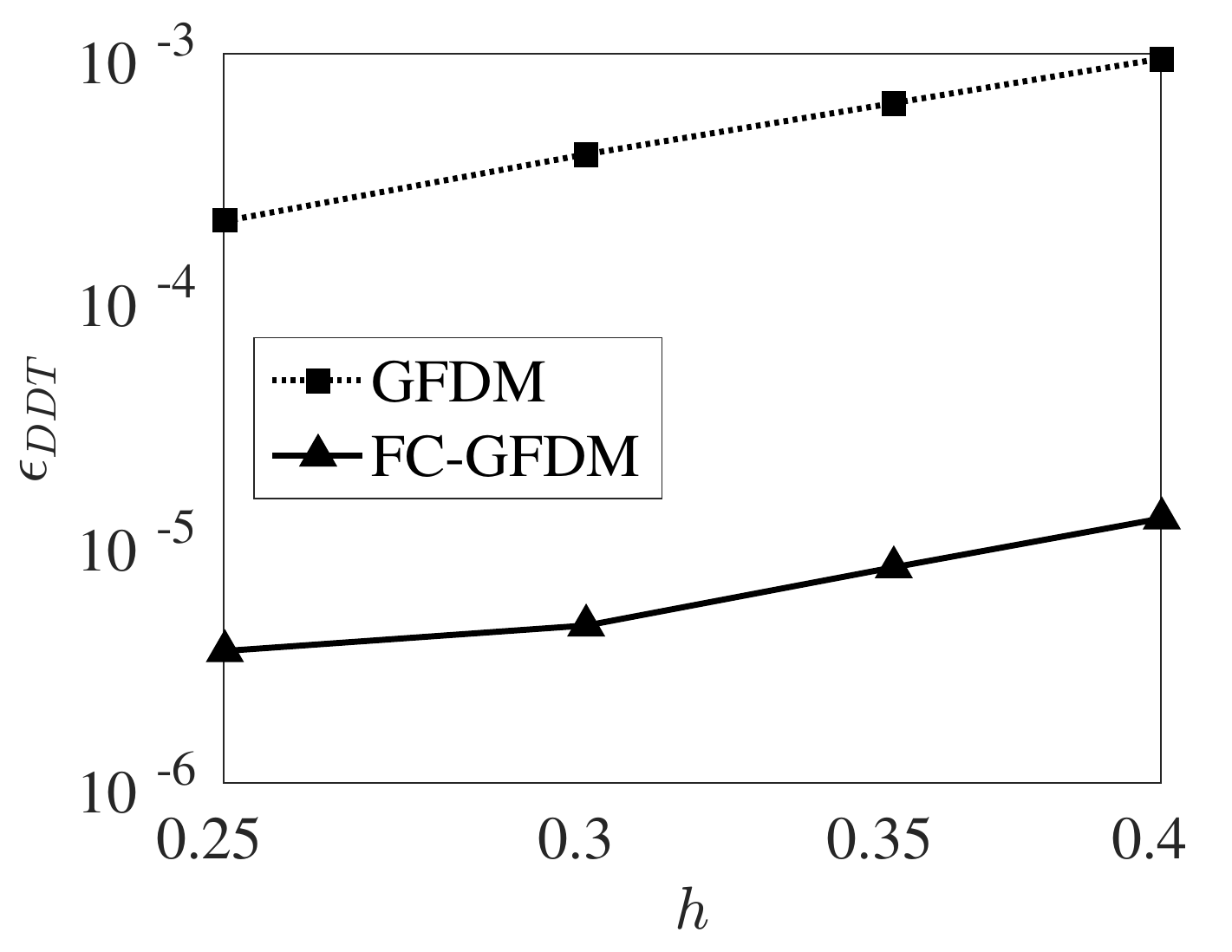}
  \includegraphics[width=0.45\textwidth]{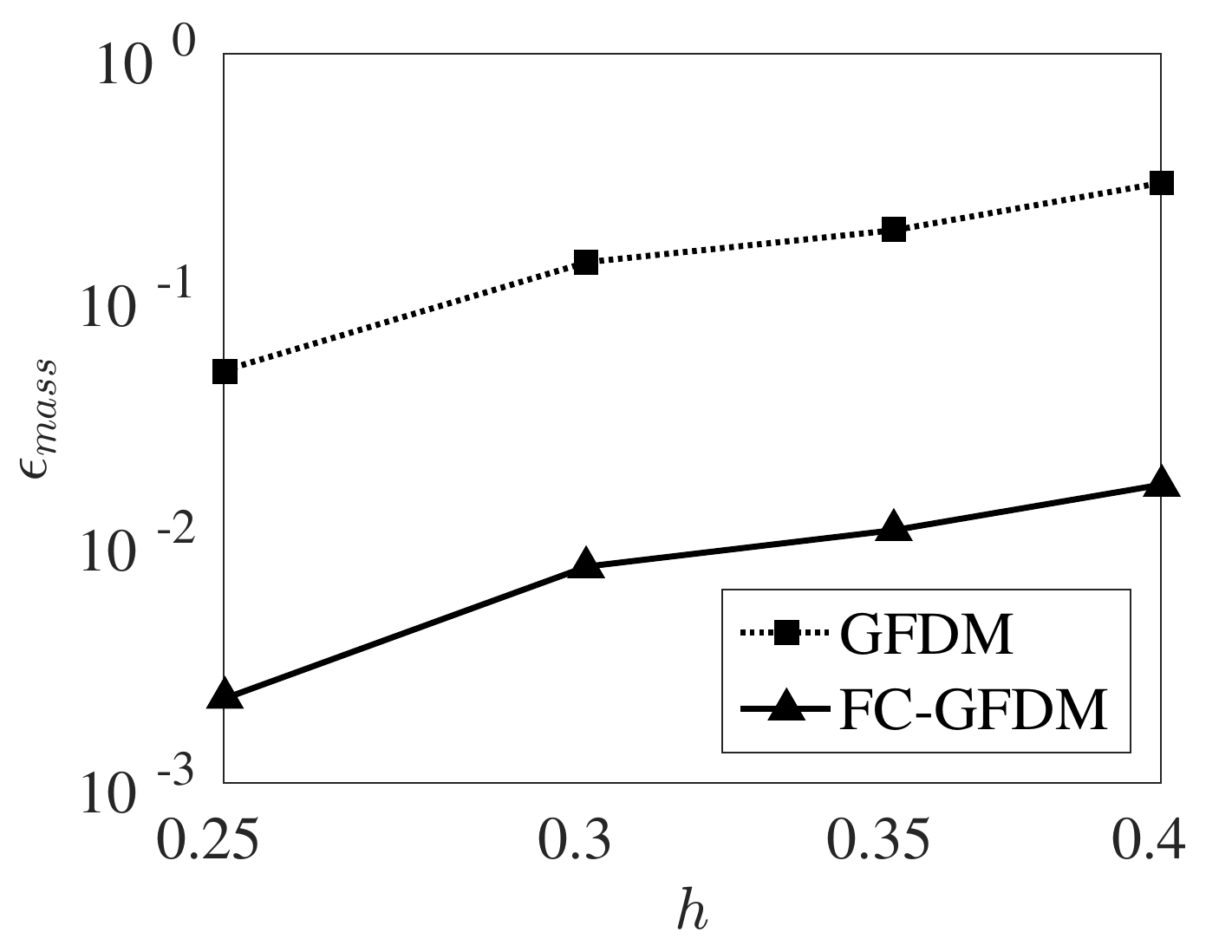}
  \caption{3D obstructed channel: convergence of solution with smoothing length $h$. Error in discrete divergence theorem~(left) and error in mass conservation~(right). The total number of points in the computational domain for the different smoothing lengths considered are shown in Table~\ref{tab:hvsN_3DChannel}.}%
  \label{Fig:3DChannel_h}
\end{figure}

The convergence of $\epsilon_{DDT}$ averaged over all time steps and $\epsilon_{mass}$ for the two methods are shown with respect to a constant time-step $\Delta t$ in Figure~\ref{Fig:3DChannel_dt}. A smoothing length of $h = 0.3$ is used, which corresponds to an initial number of points $N=20\,043$ in the entire domain. The same convergence with respect to the smoothing length $h$ is shown in Figure~\ref{Fig:3DChannel_h}. A varying time-step is used according to Eq.\,\eqref{Eq:VarDt} with $C_{\Delta t}=0.3$. Table~\ref{tab:hvsN_3DChannel} shows the total number of initial points in the entire domain for the different smoothing lengths considered. Figures \ref{Fig:3DChannel_dt} and \ref{Fig:3DChannel_h} illustrate that for a fixed $h$ and $\Delta t$, the results produced by the new FC-GFDM exhibit a much smaller error in a global divergence theorem. That, in turn, results in significantly smaller errors in mass conservation. Similar to the 2D case presented earlier, the errors in FC-GFDM are lower by an order of magnitude.
\begin{table}[!htbp]
	\caption{$h$ vs. $N$ for 3D Channel}
	\centering
	\tabsize
	\label{tab:hvsN_3DChannel}
	\begin{tabular}{ll}
	\toprule
	Smoothing Length $h$ & Total number of initial Points $N$ \\
	\midrule
	$0.4$  &  $\phantom{a}9\,707$  \\
	$0.35$ &  $13\,680$  \\
	$0.3$  &  $20\,043$  \\
	$0.25$ &  $30\,796$  \\
	\bottomrule
	\end{tabular}
\end{table}

\begin{figure}
  \centering
  \includegraphics[width=0.45\textwidth]{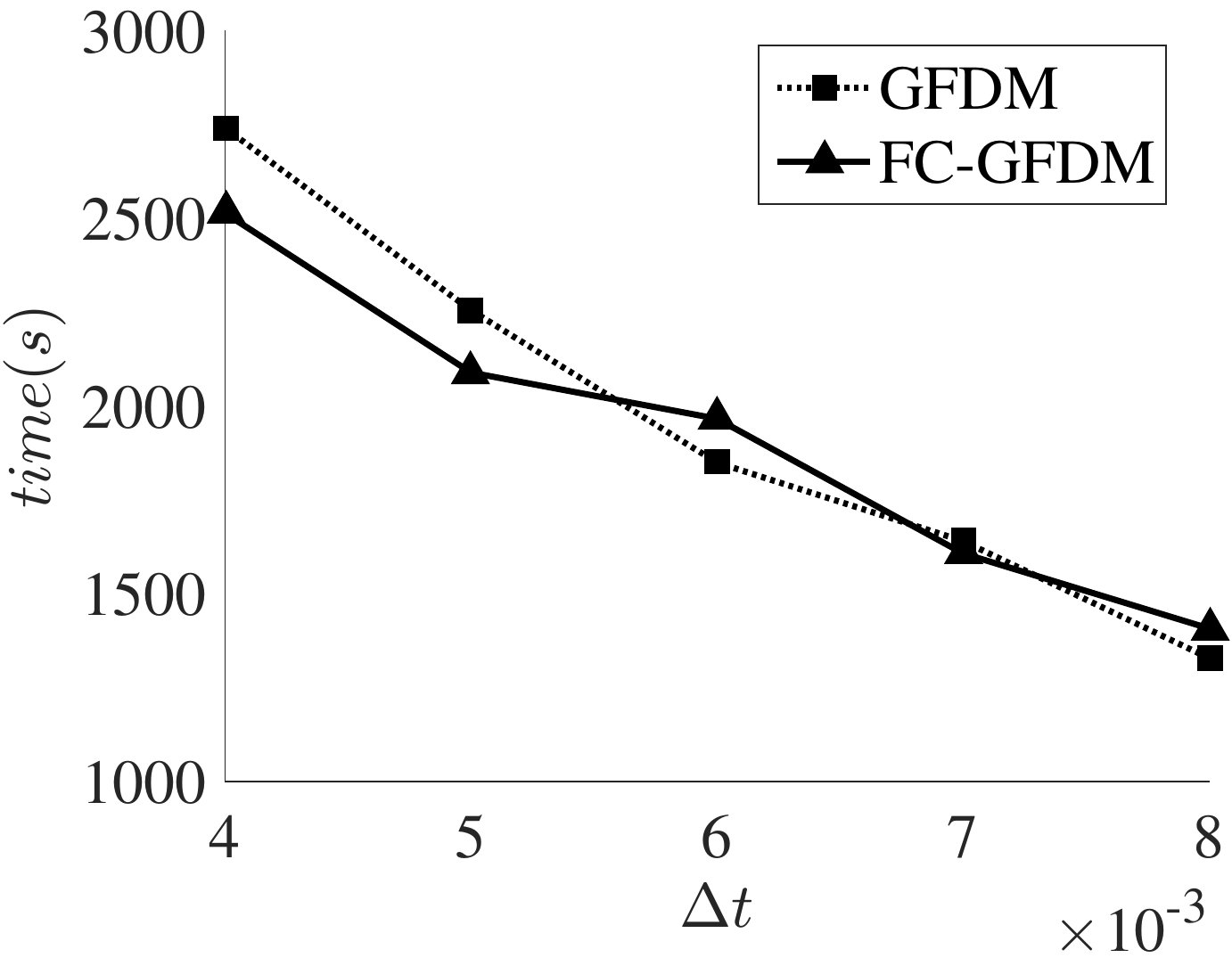}
  \includegraphics[width=0.45\textwidth]{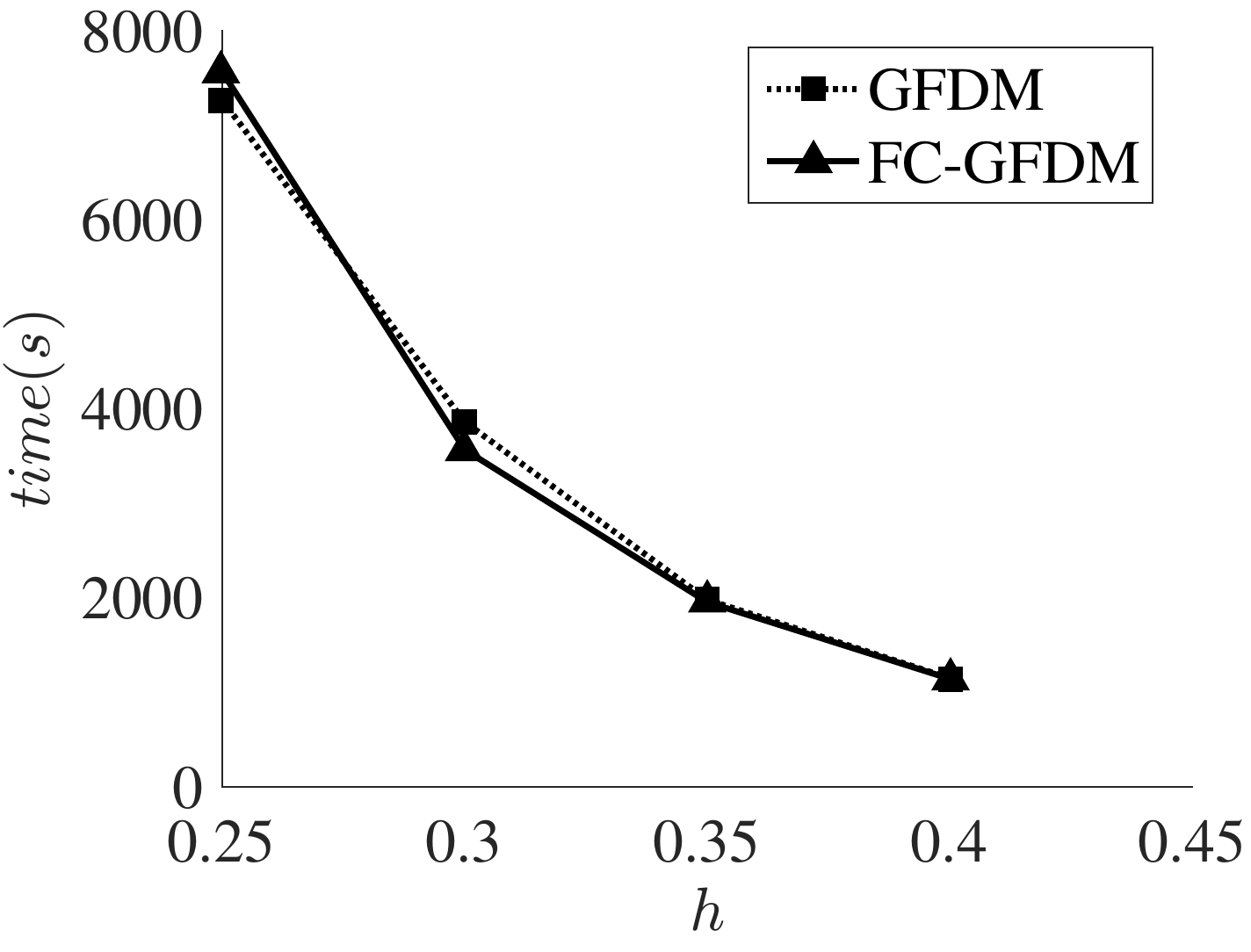}
  \caption{3D obstructed channel: clock time}
  \label{Fig:3DChannel_TimeTaken}%
\end{figure}

The time taken for the simulation for both methods are shown in Figure~\ref{Fig:3DChannel_TimeTaken}. The values of time taken shown in the figure represent the total time taken by each simulation, including the setup of the initial point cloud, the computation of differential operators and the solving of the large sparse linear systems at each time-step, and the post-processing integrations. The simulations were carried out using Fortran and were run serially on an Intel XeonE5-2670 CPU rated at 2.60GHz. Time comparisons were done under the exact same conditions. With the exception of the differential operators, both methods use the same subroutines; while the differential operators of both GFDM and FC-GFDM use the same implementation of a QR-decomposition. Figure~\ref{Fig:3DChannel_TimeTaken} illustrates that, for the same $h$ and $\Delta t$, both methods take a similar amount of time. For almost no additional computational effort, the FC-GFDM produces much smaller errors in mass conservation. Further, Figures~\ref{Fig:3DChannel_dt} - \ref{Fig:3DChannel_TimeTaken} show that the FC-GFDM takes significantly lesser time to achieve a certain tolerance of mass conservation. The addition of the flux conservation condition results in a slightly larger computation time for the differential operators of the FC-GFDM compared to those of the classical GFDM. However, as stated earlier, the effect of this on the overall simulation time is not significant. In fact, occasionally, the overall simulation time can be slightly lesser for the FC-GFDM, as illustrated in Figure~\ref{Fig:3DChannel_TimeTaken}. This can be explained by a possible faster convergence of the large sparse linear systems in the FC-GFDM case. 

\subsection{Sloshing}
\begin{figure}
  \centering
  \includegraphics[width=0.4\textwidth]{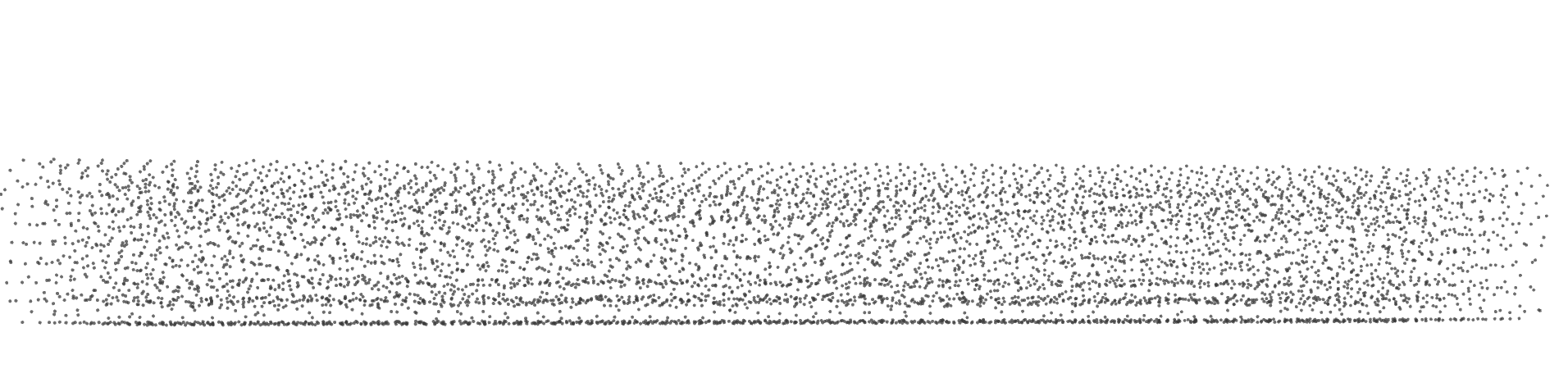}
  \includegraphics[width=0.4\textwidth]{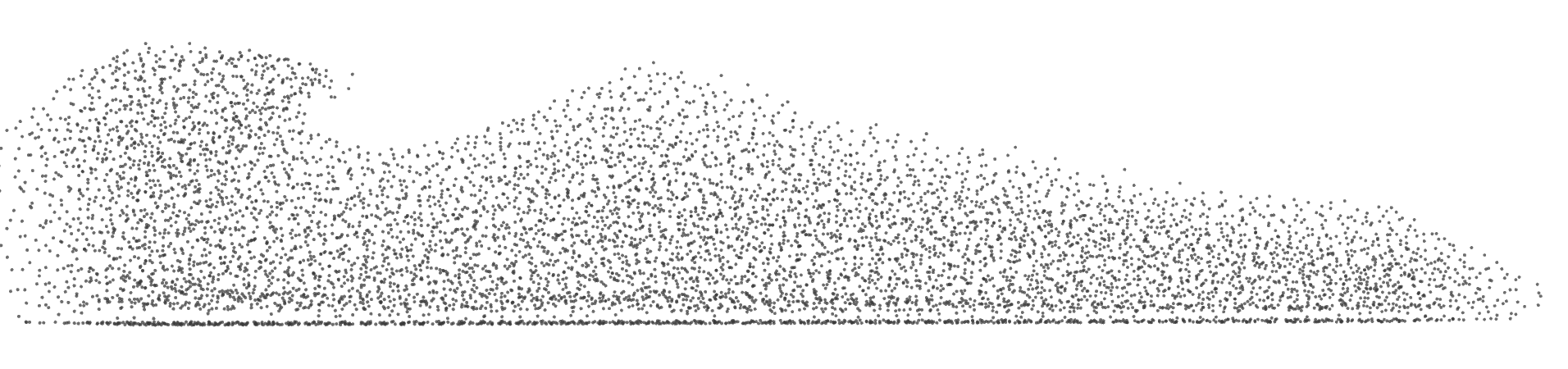}
  \vfill
  \includegraphics[width=0.4\textwidth]{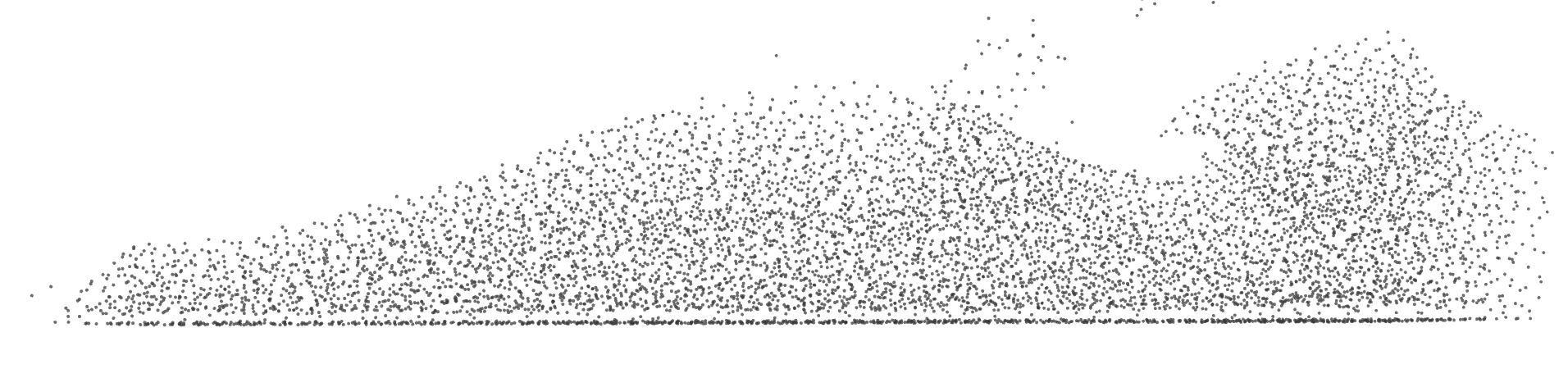}
  \includegraphics[width=0.4\textwidth]{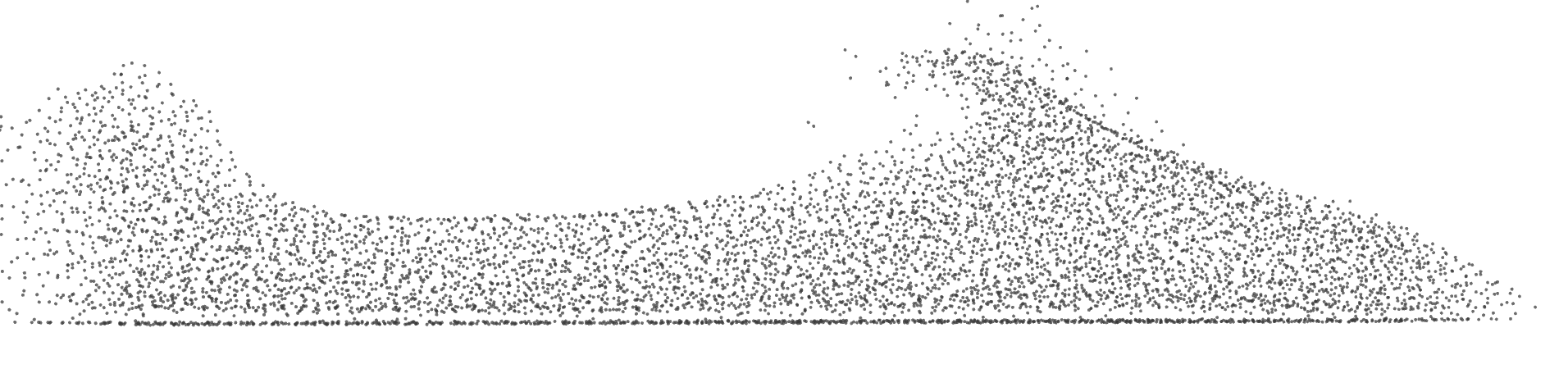}
  \caption{Sloshing. $t = 0s$~(top left), $t = 0.75s$~(top right), $t = 2.0s$~(bottom left) and $t = 3.2s$~(bottom right)}
  \label{Fig:Sloshing}%
\end{figure}

To illustrate the use of the new FC-GFDM differential operators on moving domains and problems with free surfaces, we consider the three dimensional sloshing of water in a rectangular box  as shown in Figure~\ref{Fig:Sloshing}.

\begin{figure}[!htbp]
  \centering
  \includegraphics[width=0.45\textwidth]{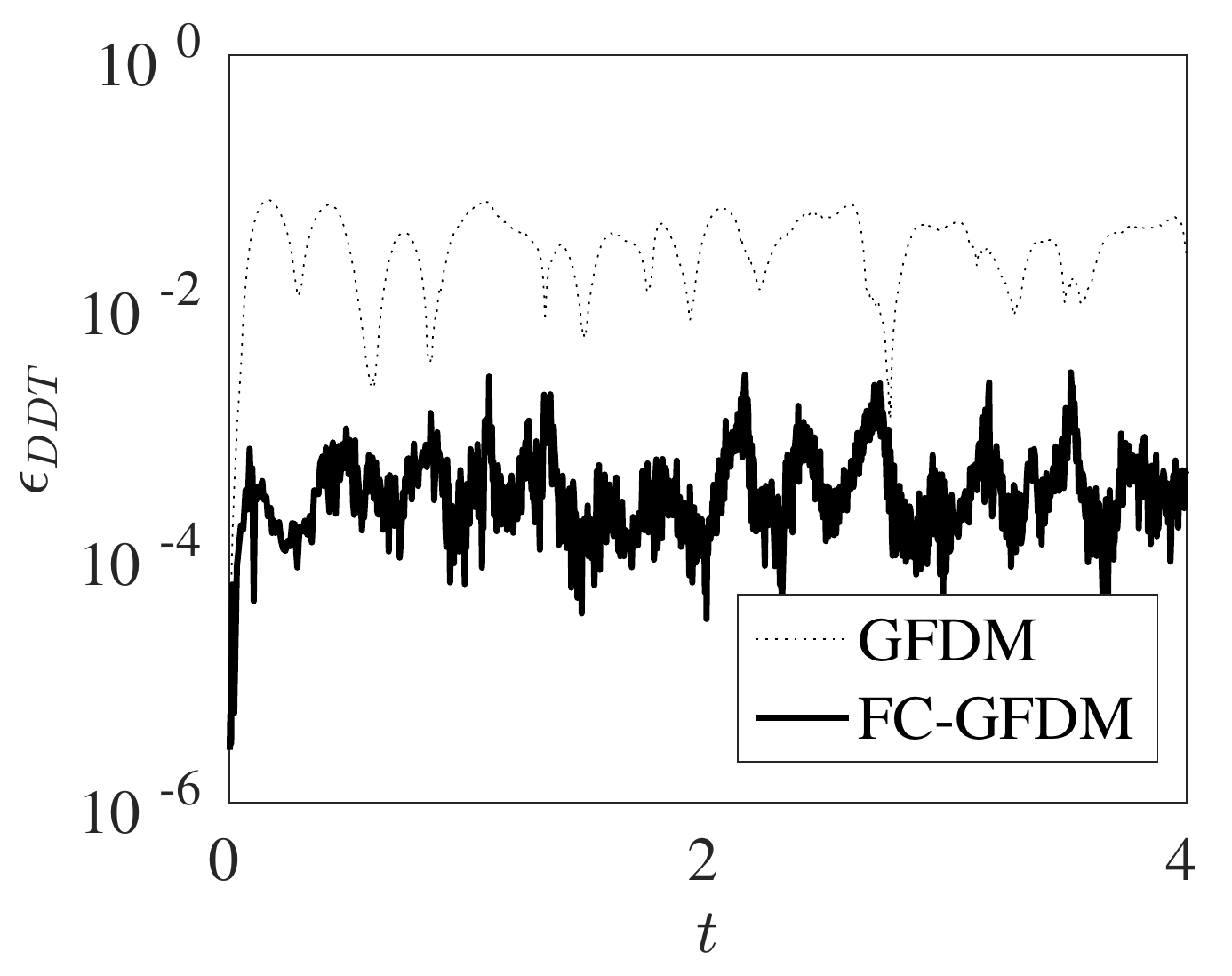}
  \includegraphics[width=0.45\textwidth]{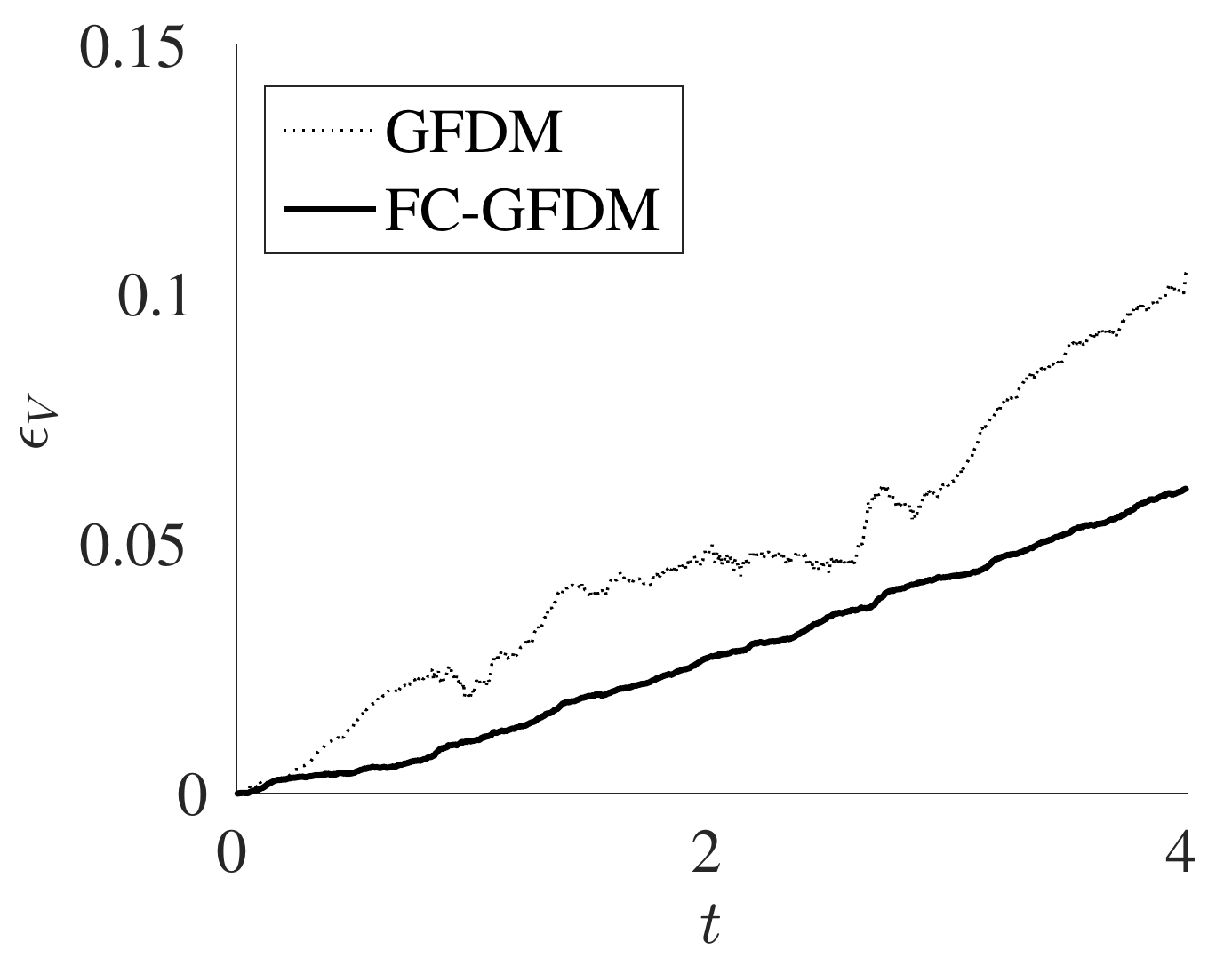}
   \caption{3D Sloshing: Time evolution of error in divergence theorem~(left) and error in mass/volume conservation~(right).}
  \label{Fig:Sloshing_TimeEvolution}
\end{figure}

The initial state is taken to be at rest. Slip boundary conditions are used at the walls for the velocity. Free surface boundary conditions are applied at the boundaries
\begin{align}
	\vec{t}_1^T\cdot\mathbf{S}\cdot\vec{n} &= 0 \,,\\
	\vec{t}_2^T\cdot\mathbf{S}\cdot\vec{n} &= 0 \,, \\
	\vec{n}^T\cdot\mathbf{S}\cdot\vec{n} &= p - p_0 \,,
\end{align}
where $p_0$ is the surrounding pressure, and $\vec{t}_1$ and $\vec{t}_2$ represent the tangential directions. Homogeneous Neumann boundary conditions are used for the pressure at the walls. A varying time-step is used according to Eq.\,\eqref{Eq:VarDt} with $C_{\Delta t}=0.3$. The movement of the box is represented in the gravitational and body forces term by setting $\vec{g}=(4\cos(7t),-10,0)$. The simulation parameters are set as $\rho = 10^3 kg/m^3$, $\eta = 10^{-3} Pa\,s$, and $t_{end} = 4s$. The error in mass conservation is measured by the change in total volume occupied by all points, since the density $\rho$ is fixed and constant throughout the domain
\begin{equation}
	\epsilon_{V}=\frac{|\int_{\Omega_0} dV - \int_{\Omega} dV|}	{\int_{\Omega_0} dV},
\end{equation}
where $\Omega_0$ is the initial domain and $\Omega$ is the domain at the time when the error is measured. For an initial number of points $N=3584$, the evolution of the error in a global divergence theorem and an error in mass conservation is shown in Figure~\ref{Fig:Sloshing_TimeEvolution}. Similar to the earlier cases, FC-GFDM shows significantly smaller errors in a global divergence theorem. That translates to smaller errors in mass conservation which is measured using the change in total volume.

%

%
%
%
%
%
\section{Conclusion}
\label{sec:Conclusion}

We presented a novel method that combines classical moving least squares approaches to meshfree differential operators and finite-volume like flux conservation over local control cells. Implicit time-integration schemes are used to discretize the PDEs, coupled with a conservation of numerical fluxes for specific fields at the previous time-level. The locally defined control cells are easy to create automatically, and do not impose any further restrictions on the quality of the point cloud. Thus, they do not introduce the drawbacks of globally defined meshes used in mesh-based methods such as finite elements and finite volumes. 

Our simulations show that the flux conserving differential operators significantly improve conservation properties of meshfree GFDMs. For the same space and time discretization, classical GFDM and the new FC-GFDM take similar simulation times, but the FC-GFDM produces smaller conservation errors. 

This method can easily be extended to incorporate sophisticated flux functions. A drawback of the method is that each numerical field that needs to be conserved has to be considered individually, which can become cumbersome for large systems of PDEs. The problem of getting general conservation, in an efficient manner, via a true discrete divergence theorem for meshfree GFDMs remains open. 

\end{document}